\documentclass[12pt,a4paper,twoside]{amsart}

\usepackage{geometry}
\usepackage{fancyhdr}
\usepackage{txfonts}
\usepackage{bbm}
\usepackage{mathrsfs}
\usepackage{amsfonts}
\usepackage{amsmath}
\usepackage{amssymb}
\usepackage{wasysym}
\usepackage{psfrag}
\usepackage{graphics,epsfig,amsmath}  
\usepackage{comment}
\usepackage{xcolor}
\usepackage{enumerate}
\usepackage{hyperref}
\usepackage[bbgreekl]{mathbbol}

\DeclareSymbolFontAlphabet{\mathbb}{AMSb}
\DeclareSymbolFontAlphabet{\mathbbl}{bbold}

\newtheorem{theorem}{Theorem}[section]

\newtheorem{proposition}{Proposition}

\newtheorem{lemma}{Lemma}[section]

\newtheorem{remark}{Remark}[section]

\newcommand{\ZZ}{\mathbb{Z}}
\newcommand{\NN}{\mathbb{N}}
\newcommand{\RR}{\mathbb{R}}
\newcommand{\CC}{\mathbb{C}}
\newcommand{\TT}{\mathbb{T}}
\newcommand{\mc}{\mathcal}

\newcommand{\red}{\textcolor{red}}

\newcommand{\scr}{\mathscr}

\allowdisplaybreaks

\numberwithin{equation}{section}

\title[Quasi-periodic solution]{Construction of quasi-periodic solutions for  delayed perturbation differential equations}

\author[X. He]{Xiaolong He}
\address{
School of Mathematical Science\\
Fudan University\\
220 Handan Road, Shanghai, 200433, People's Republic of China}

\email{hexiaolong@hnu.edu.cn}

\author[X. Yuan]{Xiaoping Yuan}
\address{
School of Mathematical Science\\
Fudan University\\
220 Handan Road, Shanghai, 200433, People's Republic of China}
\email{xpyuan@fudan.edu.cn}

\begin{document}

\vspace{0.1in}

\begin{abstract}
We employ the Craig-Wayne-Bourgain method to construct quasi-periodic
solutions for  delayed perturbation differential equations.
Our results not only implement the existing literatures on constructing
quasi-periodic solutions for DDE by the KAM method and space splitting technique,
but also provide an example of application of multi-scale analysis method to
non-selfadjoint problem.
\end{abstract}



\maketitle

\section{Introduction}
\subsection{Background and motivation}

There exist  plenty of literatures on the periodic theory for delay differential
equations (DDEs) by various methods (see  e.g.
 \cite{Guo15,GW13,HVL93,Wu96,ZG14,ZG17}).
However, as far as we know, it seems little attentions have been paid on the
quasi-periodic theory for delay systems, besides  the bifurcation and numerical arguments
on the model problems. The main difficulty in the construction of
quasi-periodic solutions (qp-solution for short) is
 the famous small divisor problem when the hyperbolicity is absent.

The study of quasi-periodic solutions for DDEs
dates back to the 1960's. In \cite{Hal69}, Halanay studied the qp-solutions
of linear DDEs
\begin{equation}\label{general}
\dot x(t)=L(t)x_{t}+f(t),
\end{equation}
where $L(t)$ and $f(t)$ are quasi-periodic with the same frequency.
By imposing some Diophantine type conditions, Halanay gave sufficient
and necessary conditions for the existence of qp-solutions
in the particular case of $L(t)$ being constant or periodic.
Later, Halanay and Yorke proposed the open problem on the existence of
qp-solutions of \eqref{general} for $L(t)$ being quasi-periodic in their survey paper \cite{HY71}.

There seems little substantial progress on the quasi-periodic theory for DDEs until
recent years. In \cite{LdlL09}, Li and de la Llave  considered
the linear DDEs with quasi-periodic perturbation
\begin{equation}\label{LdlL equ}
\dot x(t)= A x(t)+ B x(t-\tau)+ \epsilon f(\omega t, x(t), x(t-\tau), \xi).
\end{equation}
In  light of the parameterization method,
they transformed the problem on the existence of qp-solution into finding an embedding
$K_{\epsilon, \xi}:$ $\TT^{d}\rightarrow C([-\tau,0],\RR^{n})\times \TT^{d}$ such that
\begin{equation*}
F_{\epsilon,\xi}\circ K_{\epsilon,\xi}= K_{\epsilon, \xi}\circ R_{\omega},
\end{equation*}
where $R_{\omega}(\theta)=\theta+\omega$ and $F_{\epsilon,\xi}$ is the
time-one solution operator for \eqref{LdlL equ} on the phase space $C([-\tau,0],\RR^{n})\times \TT^{d}$.
To solve the above functional equation, Li and de la Llave made fully use of the space
splitting of $C([-\tau,0],\RR^{n})$ induced from the linear equation
\begin{equation*}
\dot x(t)= A x(t)+ B x(t-\tau)
\end{equation*}
and applied the Nash-Moser iterations. At each step,
they employed the exponential trichotomy to solve the linearized equation on
the tangent bundle, and
it is on the finitely dimensional center subbundles where
the small divisor problem is overcome.

Later, Li and Yuan \cite{LY12} considered the persistence of qp-solutions for autonomous DDEs
\begin{equation}\label{LY equ}
\dot x(t)= A(\xi) x + B(\xi) x(t-\tau)+ \epsilon f(x(t), x(t-\tau), \xi).
\end{equation}
With the aid of spectrum decomposition for the associated linear delay equation,
they wrote \eqref{LY equ} as an ODE on infinitely dimensional space
$\mathcal{BC}=C([-\tau,0],\RR^{n})\oplus X_{0}$, whose invariance equation on
the center subspace is
\begin{equation*}
\left\{
\begin{aligned}
&\dot \varphi(t)= \omega(\xi)+\epsilon M_{1}(\varphi, I_{1}) \Psi_{\Lambda_{1}}(0,\xi)
f(\varphi, I, y_{t}, \xi ),\\
&\dot I_{1}(t)= \epsilon M_{2}(\varphi) \Psi_{\Lambda_{1}}(0, \xi) f(\varphi, I, y_{t}, \xi ).
\end{aligned}
\right.
\end{equation*}
By a sequence change of variables (linear in  $I_2$ and $y_{t}$, i.e. the variables in
the hyperbolic direction), the authors applied the KAM technique to obtain an integrable normal form,
which guaranteed the existence of qp-solution for \eqref{LY equ}.

More recently, Li and Shang \cite{LS18} constructed qp-solution for DDEs with an elliptic type degenerate equilibrium,
whose proof  was also based on the decomposition of the extended phase space $\mathcal{BC}$
according to the spectrum of the linear delay equation.  See
 \cite{FB76_2,FB76, HdlL16b, SB03} for
more references on the application of KAM method to DDEs.
It is worthy noticing that, in \cite{FB76_2, FB76, LdlL09,  LS18,LY12}, the phase space decomposition (according to the spectrum of autonomous linear DDE) plays an important role, making it  possible to
deal with small divisor problem on the finitely dimensional center subspace.

However, when $B=0$ in \eqref{LdlL equ} and \eqref{LY equ}, the associated linear
equations  do not involve the time delay and thus become linear ODEs.
Such kinds of equations fall into the scope of the so called delayed perturbation
differential equations, taking the form of
\begin{equation*}
\dot x(t)=g(t, x(t))+\epsilon f(x_{t}, t),
\end{equation*}
 which attract lots of attentions
  (see e.g. \cite{CPS18, CPG89, FOS18,FS09, Win74} and references therein) over the years.
Now a  very natural question is whether systems \eqref{LdlL equ} and \eqref{LY equ}
with $B=0$ still have qp-solutions.
Although the equations at hand look simpler,
we are no longer able to apply the powerful space decomposition technique to
attack the nonlinearity involving time delay. Even for the linear equation
\begin{equation*}
  \dot x(t)= A x(t)+\epsilon A'(\omega t) x(t-\tau)+\epsilon g(\omega t),\quad x\in\RR^{n},
\end{equation*}
utilizing a naive quasi-periodic change of variable on $\RR^{n}$ turns out
to be difficult to solve the above equation up to $\mathcal{O}(\epsilon^{2})$.
For that reason, we have to resort to other methods when pursuing qp-solutions.

In a series of papers \cite{Bou94_IMRS, Bou97_MRL, Bou98_Ann, Bou02_Ann}, Bourgain
developed an alternative profound method which was originally proposed by
Craig and Wayne in \cite{CW93}, in order to overcome small divisor problem and  the unavailability of the
second Melnikov condition when studying Hamiltonian PDEs. In contrast with
the KAM theory,
the CWB method (named after Craig,Wayne and Bourgain)  is more flexible
in dealing with resonant cases and finds its application in the
Anderson location theory (see \cite{Spe88} and references therein), in
spectrum theory for Schr\"odinger operator \cite{Bou05_Green},
and in the construction of periodic and quasi-periodic solutions
for ODEs and PDEs (see \cite{Bou94_IMRS, Bou97_MRL,Bou98_Ann, CW93,Kri05}).

\subsection{Main result.}
The nonlinear delayed perturbation differential equation under consideration in this
paper is
\begin{equation}\label{nonlinear}
\dot x(t)= A x(t)+\epsilon f(x(t-\tau))+\epsilon g(\omega t).
\end{equation}
 The frequency $\omega$ is considered as parameters on which
parameter excision is taken such that \eqref{nonlinear} admits qp-solution
for the admissible frequency whenever the perturbation is small enough.

To begin with, we state our basic assumptions.
\begin{enumerate}
\item[\textbf{(H1)}]  The functions $f: \RR^{2n}\rightarrow \RR^{2n}$ and $g: \TT^{d}\rightarrow \RR^{2n}$ are real analytic.
\item[\textbf{(H2)}] The constant coefficient matrix $A$ admits $n$ pair of  \emph{simple} purely imaginary eigenvalues $\pm \textbf{i}\lambda_{j}$ with
$0<\lambda_{1}<\lambda_2<\cdots<\lambda_{n}$ . The associated  eigenvectors are $\{v_{j},\bar{v}_{j}\in\CC^{2n}: 1\leq j\leq n\}$ with
 $Av_{j}=\textbf{i}\lambda_{j} v_{j}$.
\end{enumerate}

It is known that the small divisor problem arise from the resonance between the tangent frequencies (the external forcing) and the normal frequencies (the elliptic eigenvalues). For that reason, we concentrate on the case of the matrix $A$
containing only purely imaginary eigenvalues. Furthermore, under assumption
$\textbf{(H2)}$, making  linear change of variables transforms \eqref{nonlinear}
into a Schrodinger-like differential equation involving time delay
\begin{equation}\label{schro-like}
   -\textbf{i} y'(t)=\Lambda y(t)+ \epsilon
f(y(t-\tau), \bar{y}(t-\tau))+\epsilon g(\omega t),\quad \Lambda=\textrm{diag}(\lambda_{1},\cdots,\lambda_{n}).
\end{equation}
We still denote the nonlinearity and inhomogeneous terms  by $f$ and $g$ respectively,
which are analytic but no longer satisfy the reality condition.

\begin{theorem}\label{Main theorem}
Consider equation \eqref{nonlinear} on $(x,\theta)\in\RR^{2n}\times\TT^{d}$ with $\omega\in \scr{U}\subset\RR^{d}$.
Assume $\emph{\textbf{(H1)}}$ and $\emph{\textbf{(H2)}}$ hold and let $0<\eta<1$.
There exists an $\epsilon^{*}>0$ and a constant $C^{*}=C^{*}(d,n)>0$ such that for all $0<\epsilon<\epsilon^{*}$,
there is a subset $\mathscr{U}_{\infty}$ of $\mathscr{U}$
satisfying \emph{mes}$[\scr{U}_{\infty}]/\emph{\textrm{mes}}[\scr{U}]\geq 1-C^{*}\eta$ such that,
for all $\omega\in\scr{U}_{\infty}$,
there is an analytic (in time $t$) quasi-periodic  solution of
equation \eqref{nonlinear}
with frequency
$\omega$ .
\end{theorem}

Theorem \ref{Main theorem} is an immediate result of the iteration lemma \ref{iteration} and a quantitative description of $\epsilon^*$ is given in
 Remark \ref{epsilon dep}. The proof of the theorem is delayed to section
 \ref{Sect IL}.

This paper is devoted to introduce  the Craig-Wayne-Bourgain method
 to construct qp-solution
for a class of DDEs for which the phase splitting technique for functional differential equations is
not applicable. In this respect, our paper is  complementary to
\cite{LdlL09, LS18, LY12} where
KAM techniques and space splitting play a role.
Besides, the method in this paper might also find its application in
population dynamics when the interactive populations under consideration
live in a fluctuating (especially quasi-periodic) environment
(see \cite{CMY10, FOS18, Zha17}).

As we shall see, one essential part of the Craig-Wayne-Bourgain method is to apply the multi-scale analysis to construct
the inverse of the linearized operator $T=D+\epsilon S$ from Newton iteration, where $D$ is a diagonal matrix.
A basic requirement for multi-scale analysis is that $S$ should be a Toeplitz matrix
with entries decaying rapidly off the diagonal.
Usually, this is indeed the case since the entries of $S$ originate from the Fourier
coefficients of some smooth function. However, due to the existence of time delay,
the Toepliz property is not at hand immediately after taking eigenvector-Fourier expansion. Nevertheless, this can be resolved by multiplying an invertible diagonal
matrix, but making the new matrix $T$ by no means of  self-adjoint.
Fortunately, the multi-scale analysis method is rather robust and independent of the
self-adjointness.

To avoid too much technique complexity, we impose the delayed perturbation
equation as simple as possible.
Some extensions can be easily obtained.
For instance, it is flexible to consider autonomous delayed perturbation
differential equation, to which it suffices to combine Lyapunov-Schimdt method
together with the multi-scale analysis in this paper.
Another extension is to waive the simplicity assumption on the purely
imaginary eigenvalues. Indeed, it is the great advantage of CWB method
to copy with problems with multiple resonances.
Moreover, we are also able to cope with equations with multiple time lags.

The paper is organized as follows. In section \ref{overview},
 an overview of the analysis,
 on the construction of qp-solutions for \eqref{schro-like}, is
 illustrated. Some technical
 preliminaries and notations are summarized at the beginning. In section
 \ref{constr. II}, we show how to construct and control the
 inverse of the linearized operator in the Newton equation, which exhibits
 the main idea of the multi-scale analysis method. In section \ref{Sect IL},
 we state and prove the iteration lemma, based on which we give
 a proof of of our main theorem.

\section{Preliminary}\label{overview}

In this section, we give an overview of  analysis on the Schrodinger-like
equation \eqref{schro-like} and prepare some useful lemmas.

Making the Ansatz that \eqref{schro-like}
does have a respond quasi-periodic solution, we transform \eqref{schro-like}
into a lattice algebraic equation by taking eigenvector-Fourier expansion.
To solve the nonlinear lattice equation, we apply Nash-Moser iterations
to effectively improve the corrections when an approximate solution is given.
To solve the Newton equation,
we need to construct the inverse of a matrix of large scale,
in which the multi-scale analysis method play a role.
As mentioned earlier, this requires
the linearized operator enjoying the Toeplitz property.


\subsection{Reduction to nonlinear lattice problem.}\label{reduction}
Making the Ansatz that equation \eqref{schro-like} does have a quasi-periodic solution
with frequency $\omega$,
 we obtain from eigenvector-Fourier expansion
 that the coefficients must satisfy the nonlinear equation
\begin{equation}\label{NL lattice}
\left\{\begin{aligned}
 (-\langle k,\omega\rangle+\lambda_{j}) \widehat{y}_{j}(k)&+
\epsilon e^{-\textbf{i}\langle k,\omega\rangle \tau}
[f_{j}(y,\bar{y})]^{\string^}(k)+\epsilon \widehat{g}_{j}(k)=0,\\
 (\langle k,\omega\rangle+\lambda_{j}) \widehat{\bar{y}}_{j}(k)
&+ \epsilon e^{-\textbf{i}\langle k,\omega\rangle \tau}
 [\overline{f_{j}(y,\bar{y})}]^{\string^}(k)+\epsilon
\widehat{\bar{g}}_{j}(k)=0,
\end{aligned}
\right.
\end{equation}
where the symbol $\string^$ describes the Fourier coefficients of
functions and the conjugation of \eqref{schro-like} is added.

Let
\begin{equation*}
  \mathcal{L}=\{m=(\mu,j,k): \mu=\pm 1, 1\leq j\leq n, k\in\ZZ^{d}\}\subset\ZZ^{d+2}.
\end{equation*}
Then equation \eqref{NL lattice} becomes a nonlinear lattice problem on $\mathcal{L}$.
To ensure the linearized operator satisfies Toeplitz property, we multiply both
sides of equations in \eqref{NL lattice} by $e^{\textbf{i}\langle k,\omega\rangle\tau}$ and
obtain
\begin{equation}\label{NL lattice+}
\mathscr{F}[y]\equiv D y+\epsilon \mathscr{W}[{y}]
+\epsilon g=0,
\end{equation}
where, with some abuse of notations,
\begin{equation*}
\begin{aligned}
  &{y}(-1,j,k)= \widehat{y}_{j}(k),\quad {g}(-1,j,k)= e^{\textbf{i}\langle k,\omega\rangle \tau}~ \widehat{g}_{j}(k),
  \\
  &{y}(+1,j,k)=\widehat{\bar{y}}_{j}(k)
  ,\quad
  {g}(+1,j,k)= e^{\textbf{i}\langle k,\omega\rangle \tau}~
  \widehat{\bar{g}}_{j}(k),
  \end{aligned}
\end{equation*}
${D}$ is a diagonal matrix with
\begin{equation*}
  {D}(m)= (\mu\langle k,\omega\rangle+\lambda_{j}) e^{\textbf{i} \langle
  k ,\omega\rangle \tau},\quad m=(\mu,j,k),
\end{equation*}
and
\begin{equation}\label{W}
\begin{aligned}
&  \mathscr{W}[{y}] (-1,j,k)=  [f_{j}(y,\bar{y})]^{\string^}(k),\\
&   \mathscr{W}[{y}] (+1,j,k)=  [\overline{f_{j}(y,\bar{y})}]^{\string^}(k).
  \end{aligned}
\end{equation}
Note that the diagonal matrix $D$ depends also on the frequency parameter
$\omega\in\scr{U}\subset\mathbb{R}^{d}$,
and the elements on the diagonal
(except $k=0$)
are no longer of real valued due to the presence of time delay.

For later application, we introduce some notations and phrases here.
The measure of a set $\scr{V}\subset\RR^{d}$, denoted by $\textrm{mes}[\scr{V}]$, always refers to the Lebesgue measure. The sharp symbol $\#$ represents the total
number of the elements for a finite set.
For any subset  $\Lambda$ of $\ZZ^{d}$, we denote
 $\mathcal{L}_{\Lambda}=\{(\mu,j,k)\in\mathcal{L}: k\in\Lambda\}$ and
write the restriction of $T$ on $\mathcal{L}_{\Lambda}$ by $T_{\Lambda}$. Given an
integer $N>0$, we denote the restriction of $T$ on $\{(\mu,j,k)\in\mathcal{L}: |k|\leq N\}$
by $T_{N}$ for short .
For a vector $y:\mathcal{L}\rightarrow \mathbb{C}$, we define the truncation
operator by
\begin{equation*}
(\Gamma_{N} y) (m)=\left\{
\begin{aligned}
& y(m),\quad |k|\leq N;\\
& 0,\quad \textrm{otherwise}.
\end{aligned}
\right.
\end{equation*}
For any $k\in\ZZ^{d}$ and any set $\Lambda\subset\ZZ^{d}$ containing $k$, we call
 $(T_{\Lambda})^{-1}$ the \emph{local inverse} of $T$ at
 $k$ with the neighborhood $\Lambda$. Given a point $k\in\ZZ^{d}$ and a set $U\subset\ZZ^{d}$,
$k+U$ denotes the translation set
 $\{k+l: l\in U\}$.
 To avoid confusion, we use the notation $A\setminus B$ for the set theoretical difference. The symbols $\vee$ and $\wedge$ describes the maximum
 and minimum operators respectively.

\subsection{Newton equation and Multi-scale analysis}

We shall apply Nash-Moser iterations to solve the nonlinear lattice equation
\eqref{NL lattice+}.
Roughly speaking, given an approximate solution $y$, we try to improve
the error by solving the Newton equation on $\mathcal{L}$
\begin{equation}\label{Newton equ}
{T} \Delta\equiv  {D} {\Delta}+\epsilon
{S}{\Delta}=-\scr{F}[{y}],
\end{equation}
where $\Delta$ is a correction of $y$, the matrix $$S=\scr{W}'[y]$$ is
the linearized operator of $\scr{W}$ at $y$.
As a infinitely dimensional matrix, the product is defined by
\[
(S\Delta)(m)=\sum_{m'\in\mathcal{L}} S(m,m') \Delta(m').
\]

Clearly, if the diagonal part ${D}$ is uniformly bounded away from zero, then ${T}$ can be inverted by a Neumann series for sufficiently small perturbation. However, when looking into the term $\mu\langle k,\omega\rangle+\lambda_{j}$, one immediately realizes that the small divisor problem prevents the diagonal from being dominant. As a result, we call those lattice points $m=(\mu,j,k)\in\mathcal{L}$ the \emph{singular sites} if ${D}(\mu,j,k)=\mathcal{O}(\epsilon)$.

To overcome the small divisor problem, we employ the multi-scale analysis method to solve the matrix equation \eqref{Newton equ}. As usual, we take advantage of
the truncation technique and consider
$$T_{N}\Delta=\Gamma_{N}\scr{F}[y]$$
instead.
The basic idea is to construct local inverses in the neighborhood of
singular sites and then to apply the coupling lemma (see Lemma \ref{Couple L})
to paste the local inverses together.

From the definition of $\scr{W}$, it follows that
the linearized operator or the matrix
 ${S}$ enjoys the Toeplitz property with respect to
$k$. More precisely, given any $\Lambda$ of $\mathbb{Z}^{d}$, $k,k'\in \Lambda$
and $q\in\ZZ^{d}$, there is
\begin{equation*}
  S((\mu,j,k+q),(\mu',j',k'+q))= S((\mu,j,k),(\mu',j',k')).
\end{equation*}
However, this is not true for the diagonal matrix $D$.
Keeping this in mind, we realize that the construction of local inverses in some
neighborhood of $k$ can be transformed into finding the local inverse of $T$ around
$k=0$, but with some modifications on the diagonal part. To make it precise,
 we introduce an extra parameter $\sigma\in\RR$ and
define
\begin{equation}\label{T(sigma)}
{T}^{\sigma}=D^{\sigma}+\epsilon S,
\end{equation}
where
\begin{equation}\label{D(sigma)}
{D}^{\sigma}(m)=  (\mu \langle k,\omega\rangle+\mu\sigma+\lambda_{j})
e^{\textbf{i}
(\langle k,\omega\rangle+\sigma)\tau},\quad m=(\mu,j,k).
\end{equation}
It then follows that
\begin{equation}\label{translation}
{T}^{\sigma}|_{q+\Lambda}((\mu,j,k+q),(\mu',j',k'+q))
={T}^{\sigma+\langle q,\omega\rangle}|_{\Lambda}
((\mu,j,k),(\mu',j',k')).
\end{equation}
for any set $\Lambda\subset \ZZ^{d}$ containing zero.

Since $\{\langle q,\omega\rangle: q\in\ZZ^{d}\}$ is dense on
the real line,
 a discussion
of $({T}_{\Lambda}^{\sigma})^{-1}$ for the full parameter range (of $\sigma$) is also applicable to the
restriction of translated intervals. For typical $\sigma$ and some $0\in\Lambda\subset\ZZ^{d}$ of large scale, we decompose
\begin{equation*}
T_{N}^{\sigma}=
\begin{pmatrix}
T_{\Omega_{1}}^{\sigma}  &   \epsilon P\\
\epsilon Q  &   T_{\Omega_{2}}^{\sigma}
\end{pmatrix},
\end{equation*}
and assume that $(T^{\sigma}_{\Omega_{1}})^{-1}$ can be established
by  induction hypothesis. Moreover, $\Omega_{2}$ is of small size, which in particular
is a singleton due to our assumption \textbf{(H2)}. Then we can formally write
\begin{equation*}
(T_{N}^{\sigma})^{-1}=
  \begin{pmatrix}
      (T_{\Omega_{1}}^{\sigma})^{-1}+\epsilon^{2} (T_{\Omega_{1}}^{\sigma})^{-1} P
      h^{-1}Q (T_{\Omega_{1}}^{\sigma})^{-1}
        & -\epsilon (T_{\Omega_{1}}^{\sigma})^{-1} P
        h^{-1}\\
      -\epsilon h^{-1} Q
      (T^{\sigma}_{\Omega_{1}})^{-1} & h^{-1}
\end{pmatrix},
\end{equation*}
where
\begin{equation*}
h=T^{\sigma}_{\Omega_{2}}-\epsilon^{2} Q (T^{\sigma}_{\Omega_{1}})^{-1} P.	
\end{equation*}
Note that the  function $h$ depends also on the frequency
parameter $\omega$. To establish and control $h^{-1}$, it suffices to exclude some
parameters $\omega$ such that $h$ stays away from zero in a reasonable way.
Due to the simplicity of our problem, there are various methods
 to achieve it.\footnote{For instance, one can apply the standard degree theory or the implicit function
theorem to calculate the zero points of $h$.} Here we still adopt
the powerful Malgrange's preparation theorem in \cite[Lemma 8.12]{Bou98_Ann}
to replace $h$ by approximated polynomials, and then take parameter
excision for those polynomials. In this fashion,
the analysis presented in this paper can be easily generalized to the case
of non-simple purely imaginary eigenvalues.
Indeed, it is great advantage of the multi-scale analysis
method to copy with problems with multiple resonance.

\subsection{Technique lemmas.}

Firstly, we write below the resolvent identity which is frequently used in this
paper.
Let $\Lambda=\Lambda_{1}+\Lambda_{2}$ be disjoint union and $\Lambda$ be bounded. The \emph{resolvent identity} for a matrix $T$ defined on $\Lambda$ is
\begin{equation}\label{RI 1}
T_{\Lambda}^{-1}=(T_{\Lambda_{1}}^{-1}+T_{\Lambda_{2}}^{-1})-
(T_{\Lambda_{1}}^{-1}+T_{\Lambda_{2}}^{-1})(T-T_{\Lambda_{1}}-T_{\Lambda_{2}})
T_{\Lambda}^{-1},
\end{equation}
and for $m\in\Lambda_{1}$, a pointwise resolvent identity is
\begin{equation}\label{RI 2}
T_{\Lambda}^{-1}(m,m')=\left\{
\begin{aligned}
	T_{\Lambda_{1}}^{-1}(m,m')- &\sum_{m_{1}\in \Lambda_{1}, m_{2}\in\Lambda_{2}}
	T^{-1}_{\Lambda_{1}}(m,m_{1}) T(m_{1}, m_{2}) T_{\Lambda}^{-1}(m_{2},m'),\quad
	\textrm{if}\quad m'\in\Lambda_{1},\\
	-& \sum_{m_{1}\in \Lambda_{1}, m_{2}\in\Lambda_{2}}
	T^{-1}_{\Lambda_{1}}(m,m_{1}) T(m_{1}, m_{2}) T_{\Lambda}^{-1}(m_{2},m'),\quad
	\textrm{if}\quad m'\not\in\Lambda_{1},
\end{aligned}
\right.	
\end{equation}
whenever the involved  inverses exist.

Next, we refer to \cite[Lemma 8.12]{Bou98_Ann} for a quantitative  version of Malgrange's preparation theorem. Roughly speaking, for an analytic function $$h(z;\omega)=z^{d}+\sum_{1\leq j<d} a_{j}(\omega) z^{j}+ \textrm{h.o.t.},
\quad |z|<\delta,~|\omega-\omega_{*}|<\rho,$$
we can find a $d$-degree polynomial $p$ and an analytic function $Q=o(1)$ such that
$h=(1+Q) p$ holds "locally" on $|z|<\delta^{-}<\delta$ and $|\omega-\omega^*|<\rho^{-}\leq \rho$.
Moreover, the derivatives of $p$ and $Q$ with respect to $\omega$ are
well controlled.
In this paper, since the singular cluster stays bounded along the iterations (hence the degree
$d$ in the function $h$ is fixed) ,
there is also a simpler version of the preparation theorem, which  can be  found in \cite[Lemma 21.4]{Kri05}.

Finally, to construct the inverse of the linearized operator at each Newton step,
we apply the coupling lemma in \cite{Bou98_Ann} and cite it here with some
according modifications.
\begin{lemma}\label{Couple L}\emph{(see \cite[Lemma 5.3 ]{Bou98_Ann})}
  Assume $T$ satisfies the off-diagonal estimate
  \begin{equation*}
    |T(m,m')|<e^{-|k-k'|^{c}},\quad k\neq k'.
  \end{equation*}
Let $\Lambda$ be an interval in $\ZZ^{d}$ and assume $\Lambda=\cup_{\alpha}\Lambda_{\alpha}$ a covering of $\Lambda$ with intervals
$\Lambda_{\alpha}$ satisfying
\begin{enumerate}[(a)]
  \item $|T_{\Lambda_{\alpha}}^{-1}(m,m')|< B$,
  \item $|T_{\Lambda_{\alpha}}^{-1}(m,m')|< K^{-C}$ for $|k-k'|>\frac{K}{100}$,
  \item for each $k\in\Lambda$, there is a $\alpha$ such that
    \begin{equation*}
      B_{K}(k)\cap\Lambda=\{k'\in\Lambda:~ |k'-k|\leq K\}\subset\Lambda_{\alpha},
    \end{equation*}
  \item $\textrm{diam}~\Lambda_{\alpha}< C' K$ for each $\alpha$.
\end{enumerate}
If $C>C_{0}(d)$ and $B,K$ are numbers satisfying
\begin{equation*}
  \log B< \frac{1}{100} K^{c},\quad\textrm{and}\quad K> K_{0}(c,C',d),
\end{equation*}
then
\begin{align*}
  |T_{\Lambda}^{-1}(m,m')|&< 2B,\\
  \textrm{and}\qquad|T_{\Lambda}^{-1}(m,m')|&< e^{-\frac{2}{3}|k-k'|^{c}},\quad \textrm{for}\quad
  |k-k'|> (100 C'K)^{\frac{1}{1-c}}.
\end{align*}
\end{lemma}

\section{Construction of $(T_{N}^{\sigma})^{-1}$}
\label{constr. II}

Throughout the rest of the paper, we write $\lessdot$ in estimates in order to
suppress various multiplicative constants, which depend only
$d,n,\tau, \mathscr{U},\lambda_{j}$ and could be made explicit, but need not be.
The norms used below are the Euclidian norm for vectors and the induced norms
for matrices. We also write $a\sim b$ to indicate $a\lessdot b$ and $b\lessdot a$.

Let $T=D+\epsilon S$ be the linearized operator of $\scr{F}$ at some
approximate solution $y$.
Recall the definition of $T^{\sigma}$ in \eqref{T(sigma)}.
Our goal in this section is to construct polynomials
to derive and control $(T^{\sigma}_{N})^{-1}$ for some large scale $N$.

The basic assumptions in this section is given below.
\begin{enumerate}

\item [\textbf{(A1)}] Melnikov condition:
\begin{equation*}
	|\langle k,\omega\rangle\pm \lambda_{j}\pm \lambda_{j'}|\geq \frac{\gamma}{|k|^{10d}}\quad \textrm{for }\quad k\in\ZZ^{d}\setminus \{0\},~
	|k|\leq 100 N_{0},~1\leq j,j'\leq n.
\end{equation*}

\item [\textbf{(A2)}] The matrix $S$ admits the off-diagonal exponential decay
\begin{equation*}\label{off decay 2}
|\partial^{\alpha}_{\omega}S(m,m')|	\lessdot  e^{-|k-k'|^{c}},\quad \textrm{for}\quad
 \alpha=0,1.
\end{equation*}

\item [\textbf{(A3)}] There exists some open set $\mathscr{U}'\subset \mathscr{U}$ of $\omega$ such that, when $\omega\in \mathscr{U}'$, the separation property holds. More precisely,
    assume $N'\in\NN$  and $k\in\ZZ^{d}$ satisfying
\begin{align*}
	& (N')^{C_{3}}\leq N,\\
	& 4 N'<|k|< (N')^{C_{3}}.
\end{align*}
Then for
$$
\sigma_{1}-\sigma_{2}=\langle k,\omega\rangle,
$$
the matrices $(T^{\sigma_{1}}_{N'})^{-1}$ and $(T^{\sigma_{2}}_{N'})^{-1}$
do not both fail the property:
\begin{equation}\label{sepa  new}
\begin{aligned}
	\| (T^{\sigma_{i}}_{N'})^{-1}\|&\leq \Phi(N'),\\
	|T^{\sigma_{i}}_{N'}(m,m')| & \leq e^{-\frac{1}{10}\rho_{N'} |k-k'|^{c}}\quad
	\textrm{for}\quad |k-k'|\geq (\rho_{N'}^{-1}\log N')^{C_{2}},
\end{aligned}
\end{equation}
where $\rho_{N'}=(\log N'/\log N_{0})^{-\frac{1}{\log C_{3}}}$
and
$\Phi(N')=(N')^{C_{1}}$.
\end{enumerate}

\subsection{Separation of singular clusters.}

Denote $N=N_{r}$ and define inductively
\begin{equation*}
N_{s-1}^{C_{3}}\sim N_{s}, \quad 1\leq s\leq r,
\end{equation*}
and $N_{0}$ is sufficiently large.
Consequently,
$$\rho_{s}=\rho_{N_{s}}\sim 10^{-s}.$$

 Let
\begin{equation}\label{Omega}
\Omega=\{-1,1\}\times \{1,2,\cdots,n\}\times [-N,N]^{d}.	
\end{equation}
Denote  for brevity
$$k_{0}\oplus N'=(k_{0}+[-N',N']^{d})\cap [-N,N]^{d}$$
for $k_{0}\in[-N,N]^{d}$.

Consider $T^{\sigma}_{N'}$ with $N'=N_{r-1}$.
 If
$(T^{\sigma}_{k_{0}\oplus N_{r-1}})^{-1}$ satisfies \eqref{sepa new} for all
$k_{0}\in [-N,N]^{d}$, then it suffices to apply the coupling lemma
(or essentially to apply the resolvent identity) to
construct the inverse of $T^{\sigma}_{N}$ and to control
$(T^{\sigma}_{N})^{-1}$. In what follows, we always treat the
worse cases, i.e., there exits some $k_{0}^{*}$ such that
$(T^{\sigma}_{k_{0}^{*}\oplus N_{r-1}})^{-1}$ does not satisfy \eqref{sepa new}.
Let $\Lambda_{r-1}$ be the set of all $k_{0}$ such that
$(T^{\sigma}_{k_{0}\oplus N_{r-1}})^{-1}$ fails
\eqref{sepa new}.
Then it follows from assumption \textbf{(A2)}
that $\Lambda_{r-1}$ is an interval (in $\ZZ^{d}$) containing
$k_{0}^*$ and is of size at most $8 N_{r-1}$. Moreover, for
$k_{0}$ falling outside of $\Lambda_{r-1}$, there is
\begin{align*}
		\| (T^{\sigma}_{k_{0}\oplus N_{r-1}})^{-1}\|&\leq \Phi(N_{r-1}),\\
	|T^{\sigma}_{k_{0}\oplus N_{r-1}}(m,m')| & \leq e^{-\frac{1}{10} \rho_{r-1}|k-k'|^{c}}\quad
	\textrm{for}\quad |k-k'|\geq (\rho_{r-1}^{-1}\log N_{r-1} )^{C_{2}}.
\end{align*}

Consider $T^{\sigma}_{N'}$ with $N'=N_{r-2}$.
Assume, also in worse case, that there
exists $k_{0}^{*}$ but in $\Lambda_{r-1}$ such that
$(T^{\sigma}_{k_{0}^*\oplus N_{r-2}})^{-1}$
fails \eqref{sepa new}. Then we can  repeat the analysis
as before
to get $\Lambda_{r-2}$ and its associated properties.
By  continuing the  process, we eventually obtain
a sequence of decreasing intervals
\begin{equation*}
[-N,N]^{d}\supset\Lambda_{r-1}\supset \Lambda_{r-2}	
\supset \cdots\supset \Lambda_{1}\supset \Lambda_{0},
\end{equation*}
where $\Lambda_{s}$ is of size at most $8 N_{s}$ for $s\geq 1$
and $\Lambda_{0}$ is of size $8 N_{0}$. By enlarging
the  size of $\Lambda_{s}$, we can also ensure
\begin{equation}\label{enlarge}
	\Lambda_{s}\supset (\Lambda_{s-1}+[-{2}N_{s},{2}N_{s}]^{d})
	\cap [-N,N]^{d},\quad s\geq 1,
\end{equation}
which results in the size of $\Lambda_{s}$ at most
${12} N_{s}$ for $s\geq 1$. Furthermore, for $k$ lying in
 $\Lambda_{s}\setminus \Lambda_{s-1}$ with
$s\geq 1$, there is
\begin{equation}\label{decay 1}
\begin{aligned}
	\|(T^{\sigma}_{k \oplus N_{s-1}})^{-1}\|\leq &\Phi(N_{s-1}),\\
	|(T^{\sigma}_{k \oplus N_{s-1}})^{-1}(m,m')|
	< & e^{-\frac{1}{10}\rho_{s-1}|k-k'|^{c}}\quad
	\textrm{for}\quad |k-k'|\geq (\rho_{s-1}^{-1}\log N_{s-1})^{C_{2}}.
\end{aligned}
\end{equation}

For $T_{\Lambda_{0}}^{\sigma}$, using assumption \textbf{(A1)}, we have
the following proposition.

\begin{proposition}
 Assume
 $$\Omega_{2}=\{m=(\mu,j,k)\in\mathcal{L}: |D^{\sigma}(m)|<\epsilon_{1}, k\in\Lambda_{0}\}$$
 is not empty. If
 \begin{equation}\label{epsilon_1_1}
\epsilon^{\frac{1}{100}}<\epsilon_{1}<\frac{1}{100} \min_{1\leq j,j'\leq n, j\neq j'}
\left\{1,
\lambda_{j}, |\lambda_{j}-\lambda_{j'}|\right\},
\end{equation}
 and
 \begin{equation}\label{epsilon_1}
 \epsilon_{1}\lessdot \gamma N_{0}^{-10d},
\end{equation}
then $\Omega_{2}$ is a singleton.
\end{proposition}

\noindent\textbf{Proof.} Recall the definition of $D^{\sigma}$ in \eqref{D(sigma)}.
For any $m=(\mu,j,k),m'=(\mu',j',k')\in\Omega_{2}$,
we have
\begin{align*}
	|\langle k'-k,\omega\rangle+\mu'\lambda_{j'}-\mu\lambda_{j}|<2 \epsilon_{1}.
\end{align*}
If $k\neq k'$, we obtain from assumption \textbf{(A1)} and \eqref{epsilon_1} that
the left hand side of the above inequality is greater than
\begin{equation*}
  \frac{\gamma}{|k'-k|^{10d}}\geq \frac{\gamma}{(10 N_{0})^{10d}}>10 \epsilon_{1},
\end{equation*}
which leads to a contradiction. With $k=k'$ and the smallness of $\epsilon_{1}$ in
\eqref{epsilon_1_1}, we get $\mu=\mu'$ and $j=j'$.
This completes the proof.
\qed

\bigskip

Suppose $\Omega_{2}=\emptyset$, the inverse of $T_{\Lambda_{0}}^{\sigma}$ can
be well controlled by applying Neumann series. We also consider the worse case that
$\Omega_{2}\neq\emptyset$ and hence a singleton, denoted by $$\Omega_{2}=\{m_{*}\}=\{(\mu_{*},j_{*},k_{*})\}.$$
Decompose
\begin{equation}\label{Lambda0_dec}
\{-1,1\}\times \{1,\cdots, n\}\times \Lambda_{0}=\Omega'_{1}
+\Omega_{2}	.
\end{equation}
 Let
$
\Omega=\Omega_{1} +\Omega_{2}$ (cf. \eqref{Omega}),
where
\begin{equation*}
	\Omega_{1}= \bigcup_{0\leq s\leq r} \Omega_{1,s}
\end{equation*}
and
\begin{align*}
\Omega_{1,0}&= \Omega'_{1}, \\
\Omega_{1,s}&=\{-1,1\}\times \{1,\cdots,n\}\times (\Lambda_{s}\setminus \Lambda_{s-1}),\quad  1\leq s\leq r-1,\\
\Omega_{1,r}&=\{-1,1\}\times \{1,\cdots,n\}\times
([-N,N]^{d}\setminus \Lambda_{r-1}).\\
\end{align*}

\subsection{Analysis of $(T^{\sigma}_{\Omega_{1}})^{-1}$.}
\label{anal_T^sig_omega_1}
Given $m\in\Omega_{1,s}$ for some $0\leq s\leq r$,
we define
\begin{equation}\label{Omega_1}
\Omega_{1}=\Gamma_{1}+\Gamma_{2},	
\end{equation}
where
\begin{equation*}
\Gamma_{1}=
\left\{
\begin{array}{cc}
	\{-1,1\}\times \{1,\cdots,n\}\times
(k\oplus N_{s-1}),\quad & \textrm{if}\quad s\geq 1,\\
 \Omega_{1,0},\quad & \textrm{if}\quad s=0.
\end{array}
\right.
\end{equation*}
In the cases of $s\geq 2$ and $s=0$, it follows from
\eqref{enlarge} and \eqref{Lambda0_dec} respectively
that
$\Gamma_{1}\cap \Omega_{2}=\emptyset.$
When $s=1$, we further assume without loss of
generality that
\begin{equation*}
	\Lambda_{0}\supset k_{*}\oplus N_{0},
\end{equation*}
which ensures $\Gamma_{1}\cap \Omega_{2}=\emptyset$.
Indeed, this can also be done by enlarging
the size of $\Lambda_{0}$.

We summarize below a proposition on the decay property of $(T_{\Gamma_{1}})^{-1}$
with  $0\leq s\leq r$.

\begin{proposition}\label{Prop Gamma}
  Under the assumptions  \emph{\textbf{(A1)-(A3)}}, if conditions \eqref{epsilon_1_1}-\eqref{epsilon_1} and
  \begin{equation}\label{N_0}
    N_{0}^{-\frac{C_{1}}{2}}\lessdot \epsilon_{1},
  \end{equation}
  hold, we have the following properties for $m\in\Omega_{1,s}\subset \Omega_{1}$
  and $0\leq s\leq r$.
  \begin{enumerate}
    \item [(i)] For $1\leq s\leq r$, there is
    \begin{equation}\label{decay 1 new}
\begin{aligned}
\|(T^{\sigma}_{\Gamma_{1}})^{-1}\| & \leq \Phi(N_{s-1}),\\
|(T^{\sigma}_{\Gamma_{1}})^{-1}(m,m')|& < e^{-\frac{1}{10}
\rho_{s-1} |k-k'|^{c}}\quad \textrm{for}\quad |k-k'|\geq (\rho_{s-1}^{-1}\log N_{s-1})^{C_{2}}.
\end{aligned}
\end{equation}

\item [(ii)] For $s=0$, there is
\begin{equation}\label{Neumann_1_new}
\begin{aligned}
	 \|(T^{\sigma}_{\Gamma_{1}})^{-1}\|\leq & ~\Phi(N_{0}),\\
	 |(T^{\sigma}_{\Gamma_{1}})^{-1}(m,m')|< &~ e^{-|k-k'|^{c}}
	 \quad \textrm{for}\quad k\neq k',
\end{aligned}
\end{equation}
and
\begin{equation}\label{neumann_2_new}
|(T^{\sigma}_{\Gamma_{1}})^{-1}(m,m')|< \frac{4}{\gamma}
(1+|k-k_{*}| )^{10d}	.
\end{equation}

  \end{enumerate}
\end{proposition}

\smallskip

Statement $(i)$ is an immediate result of \eqref{decay 1}
and the definition of $\Gamma_{1}$. The proof of statement $(ii)$ is
based on the Neumann series, which  appears frequently in
this paper. For the moment, we show a detailed proof and omit
similar arguments afterwards.

 \bigskip

\noindent\textbf{Proof.}  For $s=0$, we have $\Gamma_{1}=\Omega_{1,0}=\Omega_{1}'$,
Then, by \eqref{Lambda0_dec} and the definition of $\Omega_{2}$,  there is
$
|D^{\sigma}(\mu,j,k)|\geq \epsilon_{1}
$
for $(\mu,j,k)\in \Gamma_{1}$
and consequently
$
\|(D^{\sigma}|_{\Gamma_{1}})^{-1}\|\leq \epsilon_{1}^{-1}.
$
Writing $(T^{\sigma}_{\Gamma_{1}})^{-1}$ into Neumann series, we obtain from
\eqref{N_0} that
$
\|(T^{\sigma}_{\Gamma_{1}})^{-1}\|\leq 2\epsilon_{1}^{-1}< N_{0}^{C_{1}/2}=\sqrt{\Phi(N_{0})}
$
as long as
$
\epsilon \|S\| \epsilon^{-1}_{1}\leq \epsilon^{\frac{99}{100}}\|S\|<\frac{1}{4}.	
$
Moreover,  for $m,m'\in\Omega_{1}$ with $k\neq k'$, we have
\begin{align*}
	(T^{\sigma}_{\Gamma_{1}})^{-1}(m,m')	=\sum_{l\geq 1} (-1)^{l}
	\epsilon^{l} [(D^{\sigma})^{-1}S]^{l}(m,m') (D^{\sigma})^{-1}(m').
\end{align*}
We see that for $l\geq 1$
\begin{align*}
| [\epsilon (D^{\sigma})^{-1} S]^{l}(m,m') |	
=& \left| \sum_{m_{1},m_{2},\cdots,m_{l-1}\in\Gamma_{1}}
[\epsilon (D^{\sigma})^{-1} S](m,m_1)
\cdots [\epsilon (D^{\sigma})^{-1} S](m_{l-1},m')\right|\\
\leq & \sum_{\cdots} \frac{\epsilon}{\epsilon_{1}} e^{-|k-k_{1}|^{c}} \cdots
\frac{\epsilon}{\epsilon_{1}} e^{-|k_{l-1}-k'|^{c}}
\leq  (\frac{\epsilon}{\epsilon_{1}})^{l} e^{-|k-k'|^{c}},
\end{align*}
and then the off-diagonal exponential decay of $(T_{\Gamma_{1}}^{\sigma})^{-1}(m,m')$
in \eqref{Neumann_1_new}
follows.

It remains to verify \eqref{neumann_2_new}.
Using \eqref{epsilon_1}, we get
\begin{align*}
	D^{\sigma}(\mu,j,k)=&|\mu(\sigma+\langle k,\omega\rangle
	+\mu\lambda_{j}) e^{\textbf{{i}}(\sigma+\langle k,\omega\rangle)\tau}|\\
	=&|\sigma+\langle k_{*},\omega\rangle
	+\mu_{*}\lambda_{j_{*}}+\langle k-k_{*},\omega\rangle
	+\mu\lambda_{j}-\mu_{*}\lambda_{j_{*}}|\\
	\geq & \frac{\gamma}{|k-k_{*}|^{10d}}-\epsilon_{1}
	\geq  \frac{\gamma}{2 |k-k_{*}|^{10d}},
\end{align*}
and consequently
\begin{equation*}
	|(D^{\sigma}(\mu,j,k))^{-1}|\leq \frac{2 |k-k_{*}|^{10d}}{\gamma }.
\end{equation*}
Then it follows from the Neumann series that
\begin{equation*}
|(T^{\sigma}_{\Omega_{1,0}})^{-1}(m,m')|< \frac{4}{\gamma}
(1+|k-k_{*}| )^{10d}	.
\end{equation*}
\qed

\smallskip

Applying the resolvent identity \eqref{RI 2} to
$(T^{\sigma}_{\Omega_{1}})^{-1}$ with respect to
the decomposition \eqref{Omega_1} yields
\begin{equation}\label{dec_new_1}
\begin{aligned}
(T^{\sigma}_{\Omega_{1}})^{-1}(m,m')	
= &(T^{\sigma}_{\Gamma_{1}})^{-1}(m,m') \delta_{m,m'}
+
\sum_{m_{1}\in\Gamma_{1}, m_{2}\in\Gamma_{2}} (T^{\sigma}_{\Gamma_{1}})^{-1}
(m,m_{1})T^{\sigma}_{\Omega_{1}}(m_{1},m_{2})
(T^{\sigma}_{\Omega_{1}})^{-1}(m_{2},m'),\\
\end{aligned}
\end{equation}
where $\delta_{m,m'}$ equals one for $m=m'$ and vanishes for the rest.
It follows from \textbf{(A2)} and Proposition \ref{Prop Gamma}
that
\begin{equation*}
\begin{aligned}
|\sum_{m_{1}\in\Gamma_{1},m_{2}\in\Gamma_{2}} 	
(T^{\sigma}_{\Gamma_{1}})^{-1}
(m,&m_{1})T^{\sigma}_{\Omega_{1}}(m_{1},m_{2})
(T^{\sigma}_{\Omega_{1}})^{-1}(m_{2},m') | \\
&\leq \epsilon \sum_{
\substack{
m_{1}\in\Gamma_{1},m_{2}\in\Gamma_{2}\\
|k-k_{1}|> (\rho_{s-1}^{-1}\log N_{s-1})^{C_{2}}
}
}
e^{-\frac{1}{10}\rho_{s-1}|k-k_{1}|^{c}-|k_{1}-k_{2}|^{c}} |(T^{\sigma}_{\Omega_{1}})^{-1}(m_{2},m') |\\
&+
\epsilon \sum_{
\substack{
m_{1}\in\Gamma_{1},m_{2}\in\Gamma_{2}\\
|k-k_{1}|\leq (\rho_{s-1}^{-1}\log N_{s-1})^{C_{2}}
}
} \Phi(N_{s-1})
e^{-|k_{1}-k_{2}|^{c}} |(T^{\sigma}_{\Omega_{1}})^{-1}(m_{2},m') |\\
&= (I)+(II).
\end{aligned}
\end{equation*}

Consider the case of  $s\geq 1$.
Recall that $\Gamma_{1}\cap \Gamma_{2}=\emptyset$ and hence
$|k-k_{2}|\geq N_{s-1}$.
Then we have
\begin{equation*}
\begin{aligned}
(I)\leq \epsilon^{\frac{9}{10}}\sum_{m_{2}\in\Gamma_{2}} e^{-\frac{1}{10}\rho_{s-1}|k-k_{2}|^{c}}& |(T^{\sigma}_{\Omega_{1}})^{-1}(m_{2},m')|\\
\leq & \frac{\sqrt{\epsilon}}{2} \max_{m_{2}
\in \Omega_{1}, |k-k_{2}|\geq N_{s-1}} e^{-\frac{1}{20}\rho_{s-1}|k-k_{2}|^{c}} |(T^{\sigma}_{\Omega_{1}})^{-1}(m_{2},m') |.
\end{aligned}
\end{equation*}
For $(II)$, since the number of $k_{1}$ in the summation is less than
$(\log N_{s-1})^{ 2 C_{2} d}$, we derive
\begin{equation*}
\begin{aligned}
(II)\leq &	~\epsilon\sum_{m_{2}\in\Gamma_{2}}
(\rho_{s-1}^{-1}\log N_{s-1})^{ 2 C_{2} d}\ \Phi(N_{s-1})
e^{-|k-k_{2}|^{c}+|k-k_{1}|^{c}}
|(T^{\sigma}_{\Omega_{1}})^{-1}(m_{2},m') |\\
\leq & ~\epsilon^{\frac{9}{10}} \sum_{m_{2}\in\Gamma_{2}} e^{-\frac{1}{10}|k-k_{2}|^{c}} |(T^{\sigma}_{\Omega_{1}})^{-1}(m_{2},m') |\\
& \times \epsilon^{\frac{1}{10}}
(\rho_{s-1}^{-1}\log N_{s-1})^{ 2 C_{2} d}\ \Phi(N_{s-1})
e^{(\log N_{s-1})^{C_{2} c}} e^{-\frac{9}{10} N_{s-1}^{c} }\\
\leq & ~\epsilon^{\frac{9}{10}} \sum_{m_{2}\in\Gamma_{2}} e^{-\frac{1}{10}|k-k_{2}|^{c}} |(T^{\sigma}_{\Omega_{1}})^{-1}(m_{2},m') |\\
\leq &~ \frac{\sqrt{\epsilon}}{2} \max_{m_{2}
\in \Omega_{1},|k_{2}-k|\geq N_{s-1}} e^{-\frac{1}{20}|k-k_{2}|^{c}} |(T^{\sigma}_{\Omega_{1}})^{-1}(m_{2},m') |.
\end{aligned}
\end{equation*}
For $s=0$,
we obtain from $\epsilon^{\frac{1}{10}} N_{0}^{C_{1}}<1$ that
\begin{equation*}
\begin{aligned}
	(II)=&~\epsilon \sum_{
\substack{
m_{1}\in\Gamma_{1},m_{2}\in\Gamma_{2}\\
k_{1}=k
}
} \Phi(N_{0})
e^{-|k_{1}-k_{2}|^{c}} |(T^{\sigma}_{\Omega_{1}})^{-1}(m_{2},m') |\\
\leq & ~\epsilon^{\frac{9}{10}} \sum_{m_{2}\in\Gamma_{2}} e^{-\frac{1}{10}|k-k_{2}|^{c}} |(T^{\sigma}_{\Omega_{1}})^{-1}(m_{2},m') | \cdot \epsilon^{\frac{1}{10}} \Phi(N_{0})\\
\leq & ~\epsilon^{\frac{9}{10}} \sum_{m_{2}\in\Gamma_{2}} e^{-\frac{1}{10}|k-k_{2}|^{c}} |(T^{\sigma}_{\Omega_{1}})^{-1}(m_{2},m') |\\
\leq & ~\frac{\sqrt{\epsilon}}{2} \max_{m_{2}
\in \Omega_{1}} e^{-\frac{1}{20}|k-k_{2}|^{c}} |(T^{\sigma}_{\Omega_{1}})^{-1}(m_{2},m') |.
\end{aligned}
\end{equation*}

All together, due to the fact that
\begin{equation*}
|k-k_{*}|> N_{s-1},\quad s\geq 1,
\end{equation*}
and $\Omega_{1}$ is finite, we conclude the following lemma.

\begin{lemma}\label{Omega_1 pre-ite}
Suppose the assumptions of Proposition \ref{Prop Gamma} and
\begin{equation}\label{epsilon N_0}
  \epsilon^{\frac{1}{10}} N_{0}^{C_{1}}<1,
\end{equation}
hold. Then, for any
 $m=(\mu,j,k)\in\Omega_{1,s}$ with $0\leq s\leq r$,
 there is some $m_{2}=(\mu_{2},j_{2},k_{2})\in\Omega_{1}$ satisfying
 $|k-k_{2}|\geq N_{s-1}$ such that
 \begin{equation}\label{Omega_1_decay_0}
 \begin{aligned}
   |(T^{\sigma}_{\Omega_{1}})^{-1}&(m,m')|\\
   \leq &
   \left\{
   \begin{aligned}
   \Phi(N_{r-1}\wedge |k-k_{*}|) +& \sqrt{\epsilon} e^{-\frac{1}{20}\rho_{s-1}|k-k_{2}|^{c}} |(T^{\sigma}_{\Omega_{1}})^{-1}(m_{2},m') |,\quad |k-k'|\leq (\rho_{s-1}^{-1}\log |k-k_{*}|)^{C_{2}},\\
   e^{-\frac{1}{10}
\rho_{s-1} |k-k'|^{c}} +& \sqrt{\epsilon} e^{-\frac{1}{20}\rho_{s-1}|k-k_{2}|^{c}} |(T^{\sigma}_{\Omega_{1}})^{-1}(m_{2},m') |,\quad |k-k'|> (\rho_{s-1}^{-1}\log |k-k_{*}|)^{C_{2}},
   \end{aligned}
   \right.
   \end{aligned}
 \end{equation}
 where $m'=(\mu',j',k')\in \Omega_{1}$.
\end{lemma}

\begin{remark}\label{refinement}
Note that
\begin{equation*}
\frac{1}{100}\rho_{s-1} N_{s-1}^{c/2}\geq \frac{N_{s-1}^{c/2}}{[\log_{N_{0}} N_{s-1}]^{1/\log C_{3}}}	>1
\end{equation*}
holds for $N_{0}$ large enough. The exponential decay term  $e^{-\frac{1}{10}
\rho_{s-1} |k-k'|^{c}}$ and  $e^{-\frac{1}{20}\rho_{s-1}|k-k_{2}|^{c}}$ in
\eqref{Omega_1_decay_0} can be replaced by
\begin{equation}\label{omega_1_decay_3}
	e^{-\frac{1}{20}\rho_{s-1}|k-k'|^{c}}e^{-|k-k'|^{\frac{c}{2}}}\quad
\textrm{and}\quad\sqrt{\epsilon} e^{-\frac{1}{25}\rho_{r-1} |k-k_{2}|^{c}}e^{-|k-k_{2}|^{\frac{c}{2}}} |(T^{\sigma}_{\Omega_{1}})^{-1}(m_{2},m') |
\end{equation}
respectively.
\end{remark}

Again applying \eqref{Omega_1_decay_0} with the improved estimate \eqref{omega_1_decay_3} to $(T^{\sigma}_{\Omega_{1}})^{-1}(m_{2},m')$, we
have
\begin{align*}
	|(T^{\sigma}_{\Omega_{1}})^{-1}(m,m')|
\leq& \Phi(N_{r-1}\wedge |k-k_{*}|)
+ \sqrt{\epsilon} e^{-|k-k_{2}|^{c/2}} \Phi(N_{r-1}\wedge |k_{2}-k_{*}|)\\
& + (\sqrt{\epsilon})^{2} e^{-|k-k_{2}|^{c/2}-|k_{2}-k_{3}|^{c/2}} |(T^{\sigma}_{\Omega_{1}})^{-1}(m_{3},m')|
\end{align*}
for some $m_{3}\in\Omega_{1}$.
Since $|k_{2}-k_{*}|\leq |k-k_{*}|+|k_{2}-k|$, there is
\begin{equation*}
\Phi(|k_{2}-k_{*}|)= |k_{2}-k_{*}|^{C_{1}}
\leq |k-k_{*}|^{2C_{1}}+ C_{4} |k_{2}-k|^{2C_{1}}
=\Phi(|k-k_{*}|^{2})+C_{4} \Phi(|k_{2}-k|)^{2},
\end{equation*}
where $C_{4}$ is a sufficiently large number depending
only on $C_{2}$.
It then follows that
\begin{equation*}
\sqrt{\epsilon} e^{-|k-k_{2}|^{c/2}} \Phi(N_{r-1}\wedge |k_{2}-k_{*}|)
\lessdot \sqrt{\epsilon}\ \Phi(|k-k_{*}|^{2}\wedge N_{r-1})
\end{equation*}
and consequently
\begin{equation*}
\begin{aligned}
	|(T^{\sigma}_{\Omega_{1}})^{-1}(m,m')|
\leq & \Phi(N_{r-1}\wedge |k-k_{*}|)+
C\sqrt{\epsilon}\ \Phi(|k-k_{*}|^{2}\wedge N_{r-1})\\
& + (\sqrt{\epsilon})^{2} e^{-|k-k_{3}|^{c/2}}
|(T^{\sigma}_{\Omega_{1}})^{-1}(m_{3},m')|.
\end{aligned}
\end{equation*}
Now it is clear the iterations can be successively proceeded
and simple induction arguments yield
\begin{equation*}
|(T^{\sigma}_{\Omega_{1}})^{-1}(m,m')|
< \Phi(|k-k_{*}|^{2}\wedge N_{r-1}) \left(1+ \sum_{j\geq 1} (
C \sqrt{\epsilon})^{j}\right)< 2 \Phi(|k-k_{*}|\wedge N_{r-1})^{2}.
\end{equation*}

If $|k-k'|>  (5 ~(\rho_{s-1})^{-1}\log |k-k_{*}|)^{C_{2}}$, we divide into
two cases.

\bigskip

\noindent\textbf{Case 1:} $|k_{2}-k'|\leq (\rho_{s_{2}-1}^{-1}\log |k_{2}-k_{*}|)^{C_{2}}$, where $k_{2}\in \Lambda_{s_{2}}\setminus \Lambda_{s_{2}-1}$.
\smallskip

Recall the definition of $\rho_{s}=\rho_{N_{s}}$ in assumption \textbf{(A3)}. We see from $|k_{2}-k_{*}|<N_{s_{2}-1}$ that
\begin{equation*}
  |k_{2}-k'|\leq (\rho_{s_{2}-1}^{-1}\log |k_{2}-k_{*}|)^{C_{2}}
  \leq (\log |k_{2}-k_{*}|)^{C_{2}(1+\frac{1}{\log C_{3}})}.
\end{equation*}

From \eqref{Omega_1_decay_0} and \eqref{omega_1_decay_3},
we have
\begin{equation*}
\begin{aligned}
	|(T^{\sigma}_{\Omega_{1}})^{-1}(m,m')|\leq &
e^{-\frac{1}{20}\rho_{r-1}|k-k'|^{c}}e^{-|k-k'|^{c/2}}+ \sqrt{\epsilon} e^{-\frac{1}{25}\rho_{r-1}|k-k_{2}|^{c}} e^{-|k-k_{2}|^{c/2}}  \Phi(N_{r-1}\wedge |k_{2}-k_{*}|)^{2}\\
& + (\sqrt{\epsilon})^{2} e^{-\frac{1}{25}\rho_{r-1}|k-k_{2}|^{c}-\frac{1}{25}\rho_{r-1}|k_{2}-k_{3}|^{c}}
e^{-|k-k_{2}|^{c/2}-|k_{2}-k_{3}|^{c/2}} |(T^{\sigma}_{\Omega_{1}})^{-1}(m_{3},m')|.
\end{aligned}
\end{equation*}

Simple computation gives
\begin{equation*}
\begin{aligned}
|k_{2}-k_{*}|\leq & |k_{2}-k'|+|k-k'|+|k-k_{*}|	\\
\leq & (\log |k_{2}-k_{*}|)^{C_{2}(1+\frac{1}{\log C_{3}})}+ |k-k'|+ 10^{\frac{1}{5}\rho_{s-1} |k-k'|^{\frac{1}{C_{2}}}}\\
\leq & \frac{1}{10} |k_{2}-k_{*}|+ (1+\epsilon')  10^{\frac{2}{5}\rho_{s-1}|k-k'|^{\frac{1}{C_{2}}}},
\end{aligned}
\end{equation*}
which implies
\begin{equation*}
|k_{2}-k_{*}|\leq 	\frac{10}{9}(1+\epsilon') 10^{\frac{2}{5}\rho_{s-1}|k-k'|^{\frac{1}{C_{2}}}}.
\end{equation*}
Consequently,
\begin{equation*}
\begin{aligned}
|k-k_{2}|\geq &|k-k'|-|k'-k_{2}|>|k-k'|-(\log |k_{2}-k_{*}|)^{C_{2}}\\
>& [1-(\frac{3}{5})^{C_{2}}]\cdot |k-k'|>\frac{3}{4} |k-k'|,
\end{aligned}
\end{equation*}
if $C_{2}$ is large enough.

Then we have
\begin{equation*}
\begin{aligned}
	 \sqrt{\epsilon}\ e^{-\frac{1}{25}\rho_{r-1}|k-k_{2}|^{c}} &e^{-|k-k_{2}|^{c/2}}\Phi(N_{r-1}\wedge |k_{2}-k_{*}|)\\
<& \sqrt{\epsilon} e^{-\frac{1}{30}\rho_{r-1}|k-k'|^{c} }
e^{-(\frac{3}{4})^{c/2} |k-k'|^{c/2}}
3 ^{C_{1}}\ 10^{C_{1}|k-k'|^{\frac{1}{C_{2}}}}
\lessdot \sqrt{\epsilon}\ e^{-\frac{1}{40}\rho_{r-1}|k-k'|^{c}}
e^{-\frac{1}{2}|k-k'|^{c/2}}
\end{aligned}
\end{equation*}
if
\begin{equation*}\label{C_2 c '}
	\frac{2}{C_{2}}<c \ll 1.
\end{equation*}
 As a result,
there is
\begin{equation*}
\begin{aligned}
	|(T^{\sigma}_{\Omega_{1}})^{-1}(m,m')|
\leq
&e^{-\frac{1}{20}\rho_{r-1}|k-k'|^{c}}e^{-|k-k'|^{c/2}}+ C\sqrt{\epsilon}\ e^{-\frac{1}{40}\rho_{r-1}|k-k'|^{c}}e^{-\frac{1}{2}|k-k'|^{c/2}}\\ +& (\sqrt{\epsilon})^{2} e^{-\frac{1}{25}\rho_{r-1}|k-k_{3}|^{c}} e^{-|k-k_{3}|^{c/2}}|(T^{\sigma}_{\Omega_{1}})^{-1}(m_{3},m')|.
\end{aligned}
\end{equation*}

\smallskip

\noindent\textbf{Case 2:} $|k_{2}-k'|> (\rho_{s_{2}-1}^{-1}\log |k_{2}-k_{*}|)^{C_{2}}$, where $k_{2}\in \Lambda_{s_{2}}\setminus \Lambda_{s_{2}-1}$.
\smallskip

From \eqref{Omega_1_decay_0} and \eqref{omega_1_decay_3}, we have
\begin{equation*}
	\begin{aligned}
	|(T^{\sigma}_{\Omega_{1}})^{-1}(m,m')|\leq &
e^{-\frac{1}{20}\rho_{r-1}|k-k'|^{c}}e^{-|k-k'|^{c/2}}+ \sqrt{\epsilon} e^{-\frac{1}{20}\rho_{r-1}|k-k_{2}|^{c}} e^{-\frac{1}{10}\rho_{r-1}|k_{2}-k'|^{c}}
e^{-|k-k_{2}|^{c/2}-|k_{2}-k'|^{c/2}}\\
& + (\sqrt{\epsilon})^{2} e^{-\frac{1}{25}\rho_{r-1}|k-k_{3}|^{c}} e^{-|k-k_{3}|^{c/2}}|(T^{\sigma}_{\Omega_{1}})^{-1}(m_{3},m')|\\
\leq &
e^{-\frac{1}{20}\rho_{r-1}|k-k'|^{c}}e^{-|k-k'|^{c/2}}+ C\sqrt{\epsilon}\ e^{-\frac{1}{40}\rho_{r-1}|k-k'|^{c}}e^{-|k-k'|^{c/2}}\\
 & +  (\sqrt{\epsilon})^{2} e^{-\frac{1}{25}\rho_{r-1}|k-k_{3}|^{c}} e^{-|k-k_{3}|^{c/2}}|(T^{\sigma}_{\Omega_{1}})^{-1}(m_{3},m')|.
\end{aligned}
\end{equation*}

Combining the two cases above, we have
\begin{equation}\label{omega_1 iteration}
\begin{aligned}
	|(T^{\sigma}_{\Omega_{1}})^{-1}(m,m')|\leq &
e^{-\frac{1}{20}\rho_{r-1}|k-k'|^{c}}e^{-|k-k'|^{c/2}}+ C\sqrt{\epsilon}\ e^{-\frac{1}{40}\rho_{r-1}|k-k'|^{c}} e^{-\frac{1}{2}|k-k'|^{c/2}} \\
& +
(\sqrt{\epsilon})^{2} e^{-\frac{1}{25}\rho_{r-1}|k-k_{3}|^{c}} e^{-|k-k_{3}|^{c/2}}|(T^{\sigma}_{\Omega_{1}})^{-1}(m_{3},m')|.
\end{aligned}
\end{equation}
It then follows from iterations of \eqref{omega_1 iteration}
that
\begin{equation*}
	|(T^{\sigma}_{\Omega_{1}})^{-1}(m,m')|\leq
	e^{-\frac{1}{40}\rho_{r-1}|k-k'|^{c}}e^{-\frac{1}{2}|k-k'|^{c/2}} \left(1+ \sum_{j\geq 1} (
C \sqrt{\epsilon})^{j}\right)<  e^{-\frac{1}{50}\rho_{r-1}|k-k'|^{c}}
e^{-\frac{1}{2}|k-k'|^{c/2}}.
\end{equation*}

Summarizing the analysis in this subsection, we conclude
\begin{equation}\label{Omega_1_z_decay_}
\begin{aligned}
	|(T^{\sigma}_{\Omega_{1}})^{-1}(m,m')|< & \Phi(|k-k_{*}|\wedge N_{r-1})^{3},\\
	|(T^{\sigma}_{\Omega_{1}})^{-1}(m,m')|< & e^{-\frac{1}{50}\rho_{r-1}|k-k'|^{c}} e^{-\frac{1}{2}|k-k'|^{c/2}},\quad \textrm{if}\quad
	|k-k'|> (5 \rho_{s-1}^{-1}\log |k-k_{*}|)^{C_{2}},
\end{aligned}
\end{equation}
and
\begin{equation}\label{Omega_1_z_norm}
\|	(T^{\sigma}_{\Omega_{1}})^{-1}\|< \Phi(N_{r-1})^{3}.
\end{equation}

Note that the derivation of \eqref{Omega_1_z_decay_} and\eqref{Omega_1_z_norm} is based on
some fixed and real parameter $(\sigma,\omega)$.
Using Neumann series, one clearly has that \eqref{Omega_1_z_decay_}
and \eqref{Omega_1_z_norm} hold (up to a constant multiplier) within the
$\frac{1}{10~N^{2}~ \Phi(N_{r-1}^{3}) }$-neighborhood of
the initial parameter choice in the complex space,
to which in the sequel we assume the parameter restricted.
With the margin in our estimates, we still assume that
\eqref{Omega_1_z_decay_} and \eqref{Omega_1_z_norm} hold
in such a complex neighborhood. More precisely, we conclude the results
in this subsection in the following proposition.

\begin{proposition}\label{Omega 1 ite}
Fix $N$ and assume the parameter $(\sigma_{*},\omega_{*})\in\mathbb{R}\times\RR^{d}$ restricted such that
\emph{\textbf{(A1)-(A3)}} hold. If the constants $N_{0}$, $C_{1}$, $C_{2}$ and $c$
are appropriately chosen such that
\begin{equation}\label{arith cond}
\begin{aligned}
  &\epsilon^{\frac{1}{100}}<\epsilon_{1}\ll 1,	\\
&\epsilon_{1}\lessdot \gamma N_{0}^{-10d},\\
&N_{0}^{-\frac{C_{1}}{2}}\lessdot \epsilon_{1},\\
&\epsilon^{\frac{1}{10}} N_{0}^{C_{1}}<1,\\
&\frac{2}{C_{2}}<c \ll 1,
\end{aligned}
\end{equation}
hold, then \eqref{Omega_1_z_decay_} and \eqref{Omega_1_z_norm} hold for all
$(\sigma,\omega)$ lying in a $\frac{1}{10~N^{2}~ \Phi(N_{r-1}^{3}) }$-neighborhood of
$(\sigma_{*},\omega_{*})$ in $\CC\times\CC^{d}$.

\end{proposition}

\subsection{Construction of polynomials}

Take block decomposition of $T^{\sigma}_{N}$ into
\begin{equation*}
T_{N}^{\sigma}=
\begin{pmatrix}
T_{\Omega_{1}}^{\sigma}  &   \epsilon P\\
\epsilon Q  &   T_{\Omega_{2}}^{\sigma}
\end{pmatrix},
\end{equation*}
and a formal inverse of $T^{\sigma}_{N}$ should take the form of
\begin{equation*}
(T_{N}^{\sigma})^{-1}=
  \begin{pmatrix}
      (T_{\Omega_{1}}^{\sigma})^{-1}+\epsilon^{2} (T_{\Omega_{1}}^{\sigma})^{-1} P
      h^{-1}Q (T_{\Omega_{1}}^{\sigma})^{-1}
        & -\epsilon (T_{\Omega_{1}}^{\sigma})^{-1} P
        h^{-1}\\
      -\epsilon h^{-1} Q
      (T^{\sigma}_{\Omega_{1}})^{-1} & h^{-1}
\end{pmatrix},
\end{equation*}
where
\begin{equation*}
h=T^{\sigma}_{\Omega_{2}}-\epsilon^{2} Q (T^{\sigma}_{\Omega_{1}})^{-1} P.	
\end{equation*}
Moreover, the norm of $(T^{\sigma}_{N})^{-1}$ admits the following
control
\begin{equation}\label{norm control_2}	
\|(T^{\sigma}_{N})^{-1}\|\leq \Phi(N_{r-1}^{3})
+ \frac{\Phi(N_{r-1}^{6})} {|h|}.
\end{equation}

Let $\sigma_{1}=\sigma+\langle k_{*},\omega\rangle+\mu_{*}\lambda_{j_{*}}$.
Replacing  $\sigma$ by $\sigma_{1}$ in $h$ leads to
\begin{equation}\label{h new}
h(\sigma_{1},\omega)=	\mu_{*}(\sigma_{1}+\mu_{*}\epsilon \grave{S}(m_{*})-\mu_{*}
\epsilon^{2} Q (T^{\sigma_{1}}_{\Omega_{1}})^{-1} P),
\end{equation}
where $T^{\sigma_{1}}_{\Omega_{1}}=T^{\sigma}_{\Omega_{1}}$ and $\grave{S}$
denotes the diagonal part of the matrix $S$.
By Proposition \ref{Omega 1 ite}, $h$ is also analytic in $\sigma_{1}$ and
$\omega$ in a complex $\frac{1}{10~N^{2}~\Phi(N_{r-1}^{3})}$-neighborhood of
the initial  parameter choice.

We would like to employ the Malgrange's preparation theorem
in \cite[Lemma 8.1.2]{Bou98_Ann} to $h(\sigma_{1},\omega)$.
It suffices to check the conditions therein.

Consider the partial derivatives of $Q(T^{\sigma}_{\Omega_{1}})^{-1} P$ with
respect to $\sigma_{1}$ and $\omega$. Observe that
\begin{align*}
\partial (T^{\sigma_{1}}_{\Omega_{1}})^{-1}=	-(T^{\sigma_{1}}_{\Omega_{1}})^{-1}\
(\partial T^{\sigma_{1}}_{\Omega_{1}})\ (T^{\sigma_{1}}_{\Omega_{1}})^{-1},\quad
\partial= \partial_{\sigma_{1}}\
\textrm{or}\ \partial_{\omega}.
\end{align*}
Moreover, because of the smallness of the imaginary parts of $\sigma$ and $\omega$, we have
\begin{align*}
|\partial_{\sigma_{1}} T^{\sigma_{1}}_{\Omega_{1}}|=|\partial_{\sigma_{1}}
D^{\sigma}|=&|e^{\textbf{i}(\sigma+\langle k,\omega\rangle)}(1+\textbf{i}(\sigma+\langle k,\omega\rangle+\mu\lambda_{j}))|\\
\leq & 2 |1+\textbf{i}(\sigma_{1}+\langle k-k_{*},\omega\rangle
+\mu\lambda_{j}-\mu_{*}\lambda_{j_*}|
\lessdot|k-k_{*}|,
\end{align*}
and similarly (together with \textbf{(A2)} )
\begin{equation}\label{T'_omega}
\begin{aligned}
|\partial_{\omega} T^{\sigma_{1}}_{\Omega_{1}}(m,m')|\lessdot (1+|k-k_*|)^{2}
e^{-|k-k'|^{c}}.
\end{aligned}
\end{equation}
Since the off-diagonal elements is independent of $\sigma_{1}$, we have
\begin{align*}
\partial_{\sigma_{1}}(Q(T^{\sigma}_{\Omega_{1}})^{-1} P)
=-Q (T^{\sigma_{1}}_{\Omega_{1}})^{-1}\
(\partial_{\sigma_{1}} T^{\sigma_{1}}_{\Omega_{1}})\ (T^{\sigma_{1}}_{\Omega_{1}})^{-1} P.
\end{align*}
However, the derivative of $Q(T^{\sigma}_{\Omega_{1}})^{-1}P$
with respect to $\omega$ is more complicated and
reads
\begin{equation*}
\partial_{\omega}(Q(T^{\sigma}_{\Omega_{1}})^{-1}P)
=(\partial_{\omega} Q) (T^{\sigma}_{\Omega_{1}})^{-1}P
- Q (T^{\sigma_{1}}_{\Omega_{1}})^{-1}\
(\partial_{\omega} T^{\sigma_{1}}_{\Omega_{1}})\ (T^{\sigma_{1}}_{\Omega_{1}})^{-1} P
+ Q (T^{\sigma}_{\Omega_{1}})^{-1} (\partial_{\omega} P)
\end{equation*}

We only analyze the complicated term $Q (T^{\sigma_{1}}_{\Omega_{1}})^{-1}\
(\partial_{\omega} T^{\sigma_{1}}_{\Omega_{1}})\ (T^{\sigma_{1}}_{\Omega_{1}})^{-1} P$, which
takes the form of
\begin{equation}\label{perturb_1}
\begin{aligned}
 \epsilon^{2} \sum_{m_{i}\in\Omega_{1},1\leq i\leq 4}
\ S(m_*, m_{1}) (T^{\sigma_{1}}_{\Omega_{1}})^{-1}(m_{1}, m_{2})
(\partial_{\omega} T^{\sigma_{1}}_{\Omega_{1}})(m_{2},m_{3})
(T^{\sigma_{1}}_{\Omega_{1}})^{-1}(m_{3},m_{4}) S(m_{4},m_{*}).
\end{aligned}
\end{equation}
Let $\Delta=|k_{1}-k_{*}|\vee |k_{3}-k_{*}|$.
If
\begin{equation}\label{complicated case}
|k_{1}-k_{2}|\leq (5~\log \Delta)^{C_{2}(1+\frac{1}{\log C_{3}})} \quad
\textrm{and}\quad
|k_{3}-k_{4}|\leq (5~\log \Delta)^{C_{2}(1+\frac{1}{\log C_{3}})},
\end{equation}
then
\begin{equation*}
\begin{aligned}
	 |\textrm{\eqref{perturb_1}}|\lessdot & ~\epsilon^{2}
\sum_{m_{i}\in\Omega_{1}}
e^{-|k_{*}-k_{1}|^{c}} \Phi(|k_{1}-k_{*}|)^{3}
(1+|k_{2}-k_{*}|)^{2} e^{-|k_{2}-k_{3}|^{c}}
\Phi(|k_{3}-k_{*}|)^{3}
e^{-|k_{4}-k_{*}|^{c}}\\
\lessdot & ~\epsilon^{2} \sum_{m_{i}\in\Omega_{1}}e^{-\frac{1}{10}\ell_{5}(c)}
\Phi(\Delta^{6})
(1+|k_{2}-k_{*}|)^{2} e^{-\frac{9}{10}\ell_{5}(c)+ |k_{1}-k_{2}|^{c}+ |k_{3}-k_{4}|^{c}},
\end{aligned}
\end{equation*}
where
\[
\ell_{5}(c)=|k_{*}-k_{1}|^{c}+|k_{1}-k_{2}|^{c}+|k_{2}-k_{3}|^{c}
+|k_{3}-k_{4}|^{c}+|k_{4}-k_{*}|^{c}.
\]
Note that
\begin{equation*}
\begin{aligned}
& \frac{9}{10}\ell_{5}(c)-|k_{1}-k_{2}|^{c}-|k_{3}-k_{4}|^{c}	\\
= & \frac{9}{10}\left\{|k_{*}-k_{1}|^{c}+\cdots+|k_{4}-k_{*}|^{c}-\left(\sqrt[c]{\frac{10}{9}}~|k_{1}-k_{2}|\right)^{c}
-\left(\sqrt[c]{\frac{10}{9}}~|k_{3}-k_{4}|\right)^{c}
\right\}\\
\geq & \frac{9}{10}
\left\{
\Delta-2 \sqrt[c]{\frac{10}{9}} (5 \log \Delta)^{C_{2}(1+\frac{1}{\log C_{3}})}
\right\}^{c}\\
\geq & \frac{1}{100} \Delta^{c}.
\end{aligned}
\end{equation*}
Therefore, if \eqref{complicated case} holds, we have
\begin{equation}\label{1/9}
|\textrm{\eqref{perturb_1}}|< \frac{1}{9}\epsilon^{\frac{3}{2}}.	
\end{equation}
The other two cases for \eqref{perturb_1} are simpler
and  leads to the same estimate as \eqref{1/9}.
Indeed, if one of the inequalities, say the first one, in \eqref{complicated case} fails,
there is
\begin{equation*}
  |k_{1}-k_{2}|> (5 \log \Delta)^{C_{2}(1+\frac{1}{\log C_{3}})}\geq (5 \rho_{s'-1}^{-1}\log |k_{1}-k_{*}|)^{C_{2}}
\end{equation*}
and therefore the off diagonal estimate in \eqref{Omega_1_z_decay_}
can be applied.

All together, we have
\begin{equation}
	|Q (T^{\sigma_{1}}_{\Omega_{1}})^{-1}\
(\partial_{\omega} T^{\sigma_{1}}_{\Omega_{1}})\ (T^{\sigma_{1}}_{\Omega_{1}})^{-1} P|
<\frac{1}{3}\epsilon^{\frac{3}{2}},
\end{equation}

Moreover, using assumption \textbf{(A2)} and repeating the analysis (from \eqref{perturb_1} to \eqref{1/9})
to $$(\partial_{\omega} Q) (T^{\sigma}_{\Omega_{1}})^{-1}P \quad
\textrm{and} \quad Q (T^{\sigma}_{\Omega_{1}})^{-1} (\partial_{\omega} P),$$
we are able to show
\begin{equation}
	|(\partial_{\omega} Q) (T^{\sigma}_{\Omega_{1}})^{-1}P|
\vee |Q (T^{\sigma}_{\Omega_{1}})^{-1} (\partial_{\omega} P)|
< \frac{1}{3}\epsilon^{\frac{3}{2}},
\end{equation}
which further implies
\begin{equation}
	|\partial_{\omega}(Q(T^{\sigma}_{\Omega_{1}})^{-1}P)|
< \epsilon^{\frac{3}{2}}.
\end{equation}

Then the perturbation $\varphi$, defined by
\begin{equation*}
\varphi(\sigma_{1}; \omega)=\mu_{*} \epsilon ( \grave{S}(m_{*})-
\epsilon Q(T^{\sigma}_{\Omega_{1}})^{-1}P),
\end{equation*}
satisfies
\begin{equation*}
	|\partial^{\alpha} \varphi(\sigma_{1};\omega)|\leq \epsilon^{\frac{2}{3}},\quad \textrm{for}\quad \alpha\in\mathbb{N}^{d+1}, |\alpha|\leq 1.
\end{equation*}

Applying Malgrange's preparation theorem to  function $h$ we derive a
first order polynomial
$$
\tilde{p}(\sigma_{1})=\sigma_{1}+a_{0}(\omega),
$$
and a function $\tilde{q}=o(1)$ such that
\begin{equation}\label{h}
h(\sigma_{1};\omega)=\mu_{*} \tilde{p}(\sigma_{1})(1+\tilde{q}(\sigma_{1},\omega)).	
\end{equation}
Once the approximate polynomial $\tilde{p}$ is obtained, we denote the
modified polynomial
$$p(\sigma_{1})=\sigma_{1}+\textbf{Re}(a_{0}(\omega)).$$

Since the preparation theorem can only be applied locally, we eventually
obtain finitely many polynomials $\tilde{p}$ and we denote by $\mathcal{P}_{N_{0}}^{(1)}$
the set of those modified  polynomials $p$. Then
 the total number of the elements in $\mathcal{P}_{N}^{(1)}$ is bounded by
\begin{equation*}
\# \mathcal{P}^{(1)}_{N}\lessdot  N
(10 N^{2}\Phi(N_{r-1})^{3})^{2}\lessdot N^{5+6\frac{C_{1}}{C_{3}}}.
\end{equation*}

Whenever parameter excision is taken on $\omega$ such that
\begin{equation}\label{cons p 0}
|p(\sigma_{1})|>\Phi(N)^{-\frac{1}{2}}
\end{equation}
for all $p\in \mathcal{P}^{(1)}_{N}$,
we then see from
 \eqref{norm control_2}, \eqref{h} and $|p(\sigma)|<|\tilde{p}(\sigma)|$
that
\begin{equation}\label{T_N^sgm}
  \|(T^{\sigma}_{N})^{-1}\|\leq  \Phi(N)
\end{equation}
provided $C_{3}>12$.

\subsection{Off diagonal exponential decay}

In this part, we establish the off diagonal exponential
decay for $(T^{\sigma}_{N})^{-1}$,
which is an immediate result of the resolvent identity.
Recall that $\Omega=\{-1,1\}\times \{1,\cdots, n\}\times[-N,N]^{d}=\Omega_{1}+\Omega_{2}$ with
$\Omega_{2}=\{m_*\}$.
Consider $|k-k'|>(\rho_{N}^{-1}\log N)^{C_{2}}$.
Applying the resolvent identity and \eqref{Omega_1_z_decay_}
leads to
\begin{equation*}
\begin{aligned}	
|(T^{\sigma}_{\Omega})^{-1}(m,m')|\leq &|(T^{\sigma}_{\Omega_{1}})^{-1}(m,m')|
+\sum_{m_{1}\in\Omega_{1},m_{2}=m_{*}}
|(T^{\sigma}_{\Omega_{1}})^{-1}(m,m_{1})|\cdot |T^{\sigma}_{\Omega}(m_{1},m_{2})|\cdot | (T^{\sigma}_{\Omega})^{-1}(m_{2},m')	|\\
\leq & e^{-\frac{1}{50}\rho_{r-1} |k-k'|^{c}}
+\left(\sum_{|k-k_{1}|\leq (5 \rho_{s-1}^{-1}\log|k-k_{*}|)^{C_{2}}}+
\sum_{|k-k_{1}|> (5 \rho_{s-1}^{-1}\log|k-k_{*}|)^{C_{2}}}
\right)(**),
\end{aligned}
\end{equation*}
since
$$|k-k'|>(\rho_{N}^{-1}\log N)^{C_{2}}\geq (5 \rho_{s-1}^{-1}\log 2 N)^{C_{2}}
\geq (5\rho_{s-1}^{-1}\log |k-k_{*}|)^{C_{2}}.$$

Assume, for instance, that $|k-k_{*}|\geq\frac{1}{2}|k-k'|$.
(Otherwise, $|k'-k_{*}|\geq\frac{1}{2}|k-k'|$ and
the analysis below is the same.)
When $|k-k_{1}|\leq (5\rho_{s-1}^{-1}\log|k-k_{*}|)^{C_{2}}$,
we see that
\begin{equation*}
	|k-k_{1}|
	\lessdot (\log |k-k_*|)^{C_{2}(1+\frac{1}{\log C_{3}})}<\frac{1}{10} |k-k_{*}|.
\end{equation*}
and then
\begin{equation*}
|k_{1}-k_{*}|\geq |k-k_{*}|-|k_{1}-k|
> \frac{9}{10}|k-k_{*}|\geq \frac{9}{20} |k-k'|.
\end{equation*}

Therefore, we get
\begin{align*}
	\sum_{|k-k_{1}|\leq (5\rho_{s-1}^{-1}\log|k-k_{*}|)^{C_{2}}}
	(**)\leq &  \sum_{|k-k_{1}|\leq (5\rho_{s-1}^{-1}\log|k-k_{*}|)^{C_{2}}} \epsilon~
\Phi(N_{r-1})^{3} e^{-\frac{2}{3}|k-k'|^{c}} \Phi(N)\\
\leq & ~\epsilon \Phi(N)^{3} e^{-\frac{1}{6} (\log N)^{C_{2}c}}
e^{-\frac{1}{2}|k-k'|^{c}}\\
< & ~\frac{1}{2} e^{-\frac{1}{2}|k-k'|^{c}}.
\end{align*}

Moreover, using \eqref{Omega_1_z_decay_}
and \eqref{C_2 c '}, we obtain
\begin{align*}
	\sum_{|k-k_{1}|> (5\rho_{s-1}^{-1}\log|k-k_{*}|)^{C_{2}}}(**)
\leq &\sum_{|k-k_{1}|> (5\rho_{s-1}^{-1}\log|k-k_{*}|)^{C_{2}}}
\epsilon e^{-\frac{1}{50}\rho_{r-1}|k-k_{1}|^{c}}
e^{-\frac{1}{2}|k-k_{1}|^{c/2}} e^{-|k_{1}-k_{*}|^{c} }\Phi(N)\\
\leq & ~ \epsilon e^{-\frac{1}{50}\rho_{r-1}|k-k_{*}|^{c}}
e^{-\frac{1}{4}|k-k_{*}|^{c/2}}\Phi(N)
\sum_{k_{1}} e^{-\frac{1}{4}|k_{1}-k|^{c/2}}\\
\leq &~ \epsilon e^{-\frac{1}{75}\rho_{r-1}|k-k'|^{c}}
e^{-\frac{1}{8}|k-k'|^{c/2}} \Phi(N) \sum_{k_{1}} e^{-\frac{1}{4}|k_{1}-k|^{c/2}}\\
< & ~  \frac{1}{2} e^{-\frac{1}{10}\rho_{r}|k-k'|^{c}} e^{-\frac{1}{10} |k-k'|^{c/2}}.
\end{align*}

All together, we have
\begin{equation}\label{T_N^sgm off}
	|(T^{\sigma}_{\Omega})^{-1}(m,m')|< e^{-\frac{1}{10} \rho_{r} |k-k'|^{c}}
 e^{-\frac{1}{10} |k-k'|^{c/2}},\quad
	|k-k'|\geq (\rho_{N}^{-1}\log N)^{C_{2}}.
\end{equation}
which meets  the estimate in \eqref{sepa  new}.

\section{Iteration Lemma and Its Proof}\label{Sect IL}

In this section, we establish an iteration lemma, which
produces a sequence of approximate solutions $y_{j}$ of
the nonlinear lattice equation \eqref{NL lattice+} and
the associated error $\scr{F}[y_{j}]$ is successively improved.
Moreover,  it is easy to see that
the approximate solutions $\{y_{j}\}$ converges rapidly, whose
limit is then a true solution of \eqref{NL lattice+}.

Choose the constants $C_{1}, \cdots, C_{6}$ and $c$ appropriately such that
\begin{align}
&M^{c}\approx 1\\
& c>\frac{2}{C_{2}},\label{cons 0}\\
  &
(1-c) C_{5}> C_{6}>1,\label{cons 1}\\
& 1-\frac{12}{C_{3}}> \frac{2(d+4)C_{5}+2(d+5)}{C_{1}},
	\label{cons 2}\\
&C_{1}-2d C_{3}-\frac{24 C_{1}}{C_{3}}>15
.\label{cons 3}
\end{align}
For instance, we can impose in order that
$M=100$, $c=10^{-3}$, $C_{2}=4\cdot 10^{3}$, $C_{3}=100$, $C_{5}=10$, $C_{6}=5$
and $C_{1}=10^{4}d$.

Recall that the parameter $\omega$ is defined on some
open set $\mathscr{U}$ in $\RR^{d}$, whose Lebesgue measure, without loss of generality, is supposed be one.
 We now state the iteration lemma. With some abuse of notation, we  also use $j$
 to indicate the iteration step.

\begin{lemma}\label{iteration}
Let $0<\eta<1$.
Consider the nonlinear lattice equation
\eqref{NL lattice+}.
Assume that at the $j$-th Newton iteration step, there exists
an approximate solution $y_{j}$ of \eqref{NL lattice+} and an open subset $\mathscr{U}_{j}\subset\mathscr{U}_{j-1}$
such that
the following statements hold.
\begin{itemize}

\item [\textbf{(S1)}$_{j}$] For any $\omega\in\mathscr{U}_{j},$ there is
\begin{align*}
&\textrm{Supp}~y_{j}\subset \{-1,1\}\times\{1,\cdots,n\}\times[-M^{j},M^{j}]^{d},\\
& |(\partial_{\omega}^{\alpha}~ y_{j})(m)|< e^{-|k|^{c}},\quad \alpha=0,1,\\
& |\Delta_{j-1}(m)|=|y_{j}(m)-y_{j-1}(m)|< e^{-|k|^{c}} e^{-\frac{1}{10} M^{jc}}.
\end{align*}
where $m=(\mu,j,k)\in\mathcal{L}$.

\item [\textbf{(S2)}$_{j}$]
For each $\omega\in\mathscr{U}_{j}$, we have
\begin{equation*}
\|\partial_{\omega}^{\alpha} \scr{F}[{y}_{j}]\|<
e^{-2 (M^{j})^{c}},\quad \alpha=0,1	.
\end{equation*}

\item [\textbf{(S3)}$_{j}$] Let $T=\scr{F}'[y_{j}]$ be the linearized operator
of $\scr{F}$ at $y=y_{j}$. Then for
 $N\leq M^{j}$ and $\omega\in\mathscr{U}_{j}$,
\begin{equation*}
\|(T_{N})^{-1}	\|< \Phi(N)=N^{C_{1}},
\end{equation*}
and
\begin{equation*}
|(T_{N})^{-1}(m,m')|< e^{-\frac{1}{2}|k-k'|^{c}},\quad\textrm{for}\quad
|k-k'|>N^{\frac{1}{C_{6}}}.	
\end{equation*}

\item [\textbf{(S4)}$_{j}$]
Define $T^{\sigma}$ as in \eqref{T(sigma)} from
 $T=\scr{F}'[y_{j}]$.
Let  $N\leq M^{j}$ and
$$4N \leq |k|\leq N^{C_{3}}.$$
Then, if $\sigma_{1}-\sigma_{2}=\langle k,\omega\rangle$,
the matrices $(T^{\sigma_{1}}_{N})^{-1}$ and $(T^{\sigma_{2}}_{N})^{-2}$ do not both fail the property:
\begin{equation}\label{sepa ind.}
\begin{aligned}
	\| (T^{\sigma_{i}}_{N})^{-1}\|&\leq \Phi(N),\\
	|T^{\sigma_{i}}_{N}(m,m')| & \leq e^{-\frac{1}{10}{\rho_{N}} |k-k'|^{c}}\quad
	\textrm{for}\quad |k-k'|\geq (\rho_{N}^{-1}\log N)^{C_{2}},
\end{aligned}
\end{equation}
where $\rho_{N}=(\log N/\log N_{0})^{-\frac{1}{\log C_{3}}}$.

\item [\textbf{(S5)}$_{j}$] The parameter set $\mathscr{U}_{j}$ satisfies the measure
estimate
\begin{equation*}
  \textrm{mes}[\mathscr{U}_{j-1}\setminus\mathscr{U}_{j}]<\frac{1}{50}\cdot
  \frac{\eta}{4^{j-1}}.
\end{equation*}
\end{itemize}

Then,
there exists a absolute constant $\epsilon^*$ such that
if $0<\epsilon<\epsilon^*$,
there is an improved approximate solution $y_{j+1}$
and open set $\mathscr{U}_{j+1}\subset \mathscr{U}_{j}$ such that
the same statements are satisfied with $j+1$ in place of $j$.
\end{lemma}

\begin{remark}\label{j_0 rmk}
As we shall see, the iteration process starts at a sufficiently large $j_{0}$
with the initial approximate solution $y_{j_{0}}=0$.
To keep the consistency of the notations, we set the expressions $y_{j_{0}-1}=0$, $\Delta_{j_{0}-1}=0$ and $\scr{U}_{j_{0}-1}=\scr{U}$. Moreover,
\emph{\textbf{(S5)}}$_{j_{0}}$ is understood as
$
  1-\textrm{mes}[\scr{U}_{j_{0}}]\leq C(d,n) \eta.
$
\end{remark}

\begin{remark}\label{epsilon dep}
The absolute constant $\epsilon^{*}$ depends only on the universal constant $d,n,\tau,\eta,\lambda_{j}, \lambda_{j}-\lambda_{j'}$ and the analytic radius
of $f$ and $g$.
\end{remark}

Now we employ the iteration lemma to prove our main result Theorem \ref{Main theorem}.
\bigskip

\noindent\textbf{Proof of Theorem \ref{Main theorem}.}
Let $\scr{U}_{\infty}=\cap_{j\geq j_{0}} \scr{U}_{j}$. We see from
Remark \ref{j_0 rmk} and \textbf{(S5)}$_{j}$ that
$$
\textrm{mes}[\scr{U}_{\infty}]>1-C^{*}\eta.
$$
Moreover, for each $\omega\in\scr{U}_{\infty}$, the sequence of
approximate solutions $\{y_{j}\}_{j\geq j_{0}}$ converges
in $\ell^{2}(\mathcal{L})$, which can be seen from
a simple deduction
\begin{align*}
 \|y_{j'}-y_{j}\|_{2}^{2}\lessdot \sum_{k\in\ZZ^{d}} |(y_{j'}-y_{j})(m)|^{2}=\sum_{k\in\ZZ^{d}}
  |\sum_{s=j}^{j'-1} \Delta_{s}(k)|^{2}
  \leq \left(\sum_{k\in\ZZ^{d}} e^{-2 |k|^{c}}\right)\cdot
  \left(\sum_{s=j}^{j'-1} e^{-\frac{1}{10} M^{sc}}\right)^{2}.
\end{align*}
Let $y_{\infty}=\lim_{j\rightarrow\infty} y_{j}$. Then there is
$\scr{F}[y_{\infty}]=0$ due to the rapid convergence of $\{y_{j}\}_{j\geq j_{0}}$.
By the reduction arguments in section \ref{reduction},
the solution $y_{\infty}$ of \eqref{NL lattice+} defines a $\scr{C}^{\infty}$,
real-valued, quasi-periodic  function $x=x(t)$
 with frequency $\omega$, which is exactly a solution of \eqref{nonlinear}.
 Using  Corollary 2  in \cite{Nus73}, we know that the solution $x=x(t)$ is analytic
 in time $t$. This completes the proof of the theorem.
 \qed
 \bigskip

The rest of this section is devoted to the proof of the iteration lemma.
For reader's convenience,
we briefly explain the main idea behind it.

In subsection \ref{IL 0}, we start the iteration process at sufficiently large
$j_{0}$ with the trivial approximate solution $y_{j_{0}}=0$.
By imposing the Melnikov condition, the diagonal matrix is then dominated
when the perturbation is small enough. As a result, the construction of
the inverse of $T_{N}^{-1}$ is a simple application of Neumann series.
Furthermore, the separation property is also benefited from the Melnikov
 condition.
In subsection \ref{IL 1}, we employ the coupling lemma to construct
the inverse of $T_{N}$ for any $N\leq M^{j+1}$, where $T=\scr{F}'[y_{j}]$.
To obtain the  control of $(T_{N})^{-1}$, we need to exclude
some parameters of $\omega$ such that the polynomials constructed stay
away from zero.
Before solving the Newton equation, we establish the separation property
for $T^{\sigma}_{N}$ with $N\leq M^{j+1}$ in subsection \ref{IL 2},
which is validated by further excluding parameters.
As a result, we get the  desired open subset $\mathscr{U}_{j+1}\subset\mathscr{U}_{j}$
and then verify \textbf{(S5)}$_{j+1}$.
However, one should bear in mind that what we have established is
the separation property for $T^{\sigma}$ depending on the $j$-th approximate solution
$y_{j}$.
In subsection \ref{IL 3}, we construct the $(j+1)$-th  approximate solution
$y_{j+1}$ and estimate the associated new error, which verifies
the induction statement \textbf{(S1)}$_{j+1}$ and
 \textbf{(S2)}$_{j+1}$. Finally, in subsection
 \ref{IL 4}, we establish induction statement \textbf{(S3)}$_{j+1}$
 and \textbf{(S4)}$_{j+1}$. Since the associated properties have already
 been proved for $T$ depending on $y_{j}$ in
 subsections \ref{IL 1} and \ref{IL 2},  it suffices to apply
 the Neumann series directly, due to the rapid decay of the corrections $\Delta_{j}$.

\subsection{Preparation step.}\label{IL 0}
Obviously, the trivial solution $y=0$ is an approximate solution of
\eqref{NL lattice+} with error $\scr{F}[y]=\mathcal{O}(\epsilon)$.
Thus, we may let the iteration start at $y_{j_{0}}=0$ with the integer $j_{0}=j_{0}(\epsilon)$
satisfying
\begin{equation}\label{j_0}
 (M^{j_{0}})^{c}\sim \log\frac{1}{\epsilon}.
\end{equation}
Observe that $\scr{F}[0]=\epsilon \scr{W}[0]$ is independent of $\omega$.
Then statement \textbf{(S1)}$_{j_{0}}$ and  \textbf{(S2)}$_{j_{0}}$ holds.

Assume the Melnikov condition
\begin{equation}\label{Mel}
	|\langle k,\omega\rangle\pm \lambda_{j}\pm \lambda_{j'}|\geq \frac{\gamma}{|k|^{10d}}\quad \textrm{for }\quad k\in\ZZ^{d}\setminus \{0\},~
	|k|\leq (100 M^{j_{0}})^{C_{3}},~1\leq j,j'\leq n.
\end{equation}
(See also \textbf{(A1)} in section \ref{constr. II}.)
Then statement \textbf{(S3)}$_{j_{0}}$ is a simple application of Neumann series
due to the dominance of the diagonal matrix. The proof is the same to that of
 $(ii)$ in Proposition \ref{Prop Gamma} with $\epsilon_{1}$ and $N_{0}$
 replaced by $\eta N^{-10d}$ and $N$
 respectively for $N\leq M^{j_{0}}$.
 It suffices to check condition \eqref{epsilon_1_1}, which follows from \eqref{j_0}
that
 \begin{equation*}
   \epsilon_{1}=\frac{\eta}{N^{10d}}\geq \frac{\eta}{M^{10 d j_{0}}}> e^{-M^{j_{0} c}}>\epsilon^{\frac{1}{100}}.
 \end{equation*}

Consider statement \textbf{(S4)}$_{j_{0}}$. Assume on the contrary that
$\|(T^{\sigma_{1}}_{N})^{-1}\|>\Phi(N)$ and $\|(T^{\sigma_{2}}_{N})^{-1}\|>\Phi(N)$.
Then there exist $m_{1}=(\mu_{1},j_{1},k_{1})$ and $m_{2}=(\mu_{2},j_{2},k_{2})$
such that
\begin{align*}
  |\sigma_{i}+\langle k_{i},\omega\rangle+\mu_{i}\lambda_{j_{i}}|<\frac{2}{\Phi(N)},
  \quad i=1,2.
\end{align*}
This leads to
\begin{equation*}
  |\sigma_{1}-\sigma_{2}+\langle k_{1}-k_{2},\omega\rangle+\mu_{1}\lambda_{j_{1}}-\mu_{2}\lambda_{j_{2}}|
  <\frac{4}{\Phi(N)},
\end{equation*}
or equivalently
\begin{equation*}
  |\langle k+k_{1}-k_{2},\omega\rangle+\mu_{1}\lambda_{j_{1}}-\mu_{2}\lambda_{j_{2}}|
  <\frac{4}{\Phi(N)}.
\end{equation*}
Note that $2N< |k+k_{1}-k_{2}|\leq 2N^{C_{3}}$. Then, by condition \eqref{Mel},
we see  from $C_{1}=100 C_{3} d$ that
\begin{equation*}
   |\langle k+k_{1}-k_{2},\omega\rangle+\mu_{1}\lambda_{j_{1}}-\mu_{2}\lambda_{j_{2}}|
   >\frac{\eta}{(2 N^{C_{3}})^{10d}}> \frac{5}{\Phi(N)},
\end{equation*}
which turns out to be a  contradiction and thus \textbf{(S4)}$_{j_{0}}$ is valid.

Denote
$$
\scr{V}_{j_{0}}=\bigcup_{\substack{\mu,\mu'=\pm 1,\\
1\leq j,j'\leq n}}\bigcup_{1\leq |k|\leq (100 M^{j_{0}})^{C_{3}}}
\left\{\omega\in\mathscr{U}:
|\langle k,\omega\rangle+\mu \lambda_{j}+\mu' \lambda_{j'}|< \frac{\eta}{|k|^{10d}}
\right\},
$$
and
\begin{equation*}
  \scr{U}_{j_{0}}=\scr{U}\setminus \scr{V}_{j_{0}}.
\end{equation*}
Obviously, $1-\textrm{mes}[\scr{U}_{j_{0}}]=\scr{V}_{j_{0}}\lessdot \eta$
and the constant here depends only on $d$ and $n$, which verifies
statement \textbf{(S5)}$_{j_{0}}$ (see Remark \ref{j_0 rmk}).

In what follows, we assume that the iteration lemma \ref{iteration}
 holds up to step $j\geq j_{0}$.

\subsection{Construction of $(T_{N})^{-1}$}\label{IL 1}

Let $T=\scr{F}'[y_{j}]$ and
assume $M^{j}< N\leq  M^{j+1}$. We shall  construct
the inverse of $T_{N}$ by the coupling lemma
established in \cite{Bou98_Ann}(see also Lemma \ref{Couple L}).
Take
\begin{equation*}
K= N^{\frac{1}{C_{5}}},\quad B=\Phi(\frac{1}{2}N), \quad C_{5}>1.
\end{equation*}
Then for any $m=(\mu,j,k)$ with $|k|\leq 10 K$, there
is
\begin{equation*}
\|(T_{10 K})^{-1}\|< \Phi(10 K)< B,
\end{equation*}
since
$$
10 K= 10 N^{\frac{1}{C_{5}}}<10 M^{\frac{j+1}{C_{5}}}< M^{j}.
$$
For any $m=(\mu,j,k)$ with $|k|>5K$, we construct the
inverse of
$$T_{k\oplus K}= T_{k+[-K,K]^{d}}=T^{\sigma=\langle k,\omega\rangle}_{K}.$$
The covering of $\{-1,1\}\times\{1,\cdots,n\}\times[-N, N]^{d}$ is now clear.

To apply the result in section \ref{constr. II} with $y=y_{j}$, we verify the conditions therein.
Indeed, choosing $N_{0}\sim M^{j_{0}}$, the  assumption \textbf{(A1)} then follows
from our construction of $\scr{U}_{j_{0}}$. The off diagonal decay property
\textbf{(A2)} is an immediate result of the exponential decay of $y_{j}$ in \textbf{(S1)}$_{j}$. The assumption \textbf{(A3)} is just our induction
assumption \textbf{(S4)}$_{j}$. The remained arithmetical conditions \eqref{arith cond} are already
established in the preparation step.
From \eqref{T_N^sgm} in the main construction section \ref{constr. II}, we have
\begin{equation*}
\|(T^{\sigma=\langle k,\omega\rangle}_{K})^{-1}\|
\leq \Phi(K)	< B.
\end{equation*}
Of course, this requires some excision of the parameter
 $\omega$ to ensure \eqref{cons p 0}, which shall be implemented later.

Note that
\begin{equation*}
\frac{K}{100}>(10 K)^{\frac{1}{C_{6}}}	
\end{equation*}
and then we have
\begin{equation*}
|(T_{10 K})^{-1}(m,m')|<e^{-\frac{1}{2}|k-k'|^{c}}<
e^{-\frac{1}{2}(\frac{1}{100})^{c} K^{c}}< K^{-C}	
\end{equation*}
for some constant $C>0$ when $|k-k'|> \frac{K}{100}$.
Moreover, for $(T^{\sigma}_{K})^{-1}$ with $\sigma=\langle k,\omega\rangle$, we see
from \eqref{T_N^sgm off} that
\begin{equation*}
|(T^{\sigma}_{K})^{-1}(m,m')|\leq e^{-\frac{1}{10}\rho_{r(K)}
|k-k'|^{c}} e^{-\frac{1}{10}|k-k'|^{c/2}} <e^{-\frac{1}{10}(\frac{K}{100})^{c/2}}
< K^{-C},
\end{equation*}
when $|k-k'|> \frac{K}{100}$.
It  remains to check
\begin{equation*}
\log B=C_{1}\log \frac{N}{2}< \frac{N^{\frac{c}{C_{5}}}}{100}=\frac{K^{c}}{100}.
\end{equation*}
Then the coupling lemma \ref{Couple L} implies
\begin{equation}\label{coup 1}
\|(T_{N})^{-1}\|< 2B=2 \Phi(\frac{1}{2}N),	
\end{equation}
and
\begin{equation*}
|(T_{N})^{-1}(m,m')|< e^{-\frac{1}{2}|k-k'|^{c}},
\quad |k-k'|> (100 C' K)^{\frac{1}{1-c}}.	
\end{equation*}
In view of \eqref{cons 1}, or equivalently,
\begin{equation*}
\frac{1}{C_{6}}>\frac{1}{(1-c) C_{5}},
\end{equation*}
we have
\begin{equation}\label{coup 2}
|(T_{N})^{-1}(m,m')|< e^{-\frac{2}{3}|k-k'|^{c}},\quad
|k-k'|>N^{\frac{1}{C_{6}}}
\end{equation}
holds for $M^{j}< N\leq M^{j+1}$.

Now we need to take parameter separation such that for all
possible $\sigma=\langle k,\omega\rangle$, the constructed polynomial
$|p(\sigma_{1})|>\Phi(K)^{-\frac{1}{2}}$ for all $p\in\mathcal{P}^{(1)}_{K}$,
where $\sigma_{1}=\sigma+\langle k_{*},\omega\rangle+\mu_{*}\lambda_{j_{*}}=
\langle k+k_{*},\omega\rangle+\mu_{*}\lambda_{j_{*}}$.

Since $|k|>5K$ and $|k_{*}|\leq K$, we have $|k+k_{*}|>4K.$
Furthermore, it follows from Malgrange's preparation theorem that
those $a_{0}(\omega)$ in $p\in\mathcal{P}^{(1)}_{K}$ stay uniformly
bounded together with their first derivatives.
Therefore,
for any fixed $\sigma_{1}$, the measure of the excluded parameter set
\begin{equation*}
\left
\{\omega: |p(\sigma_{1})|=|\langle k+k_{*},\omega\rangle+\mu_{*}\lambda_{j_{*}}+\textbf{Re}(a_{0}(\omega))|
\leq \frac{\eta}{\sqrt{\Phi(K)}}\right\}
\end{equation*}
 is
less than $\Phi(K)^{-\frac{1}{2}}$ up to a constant multiplier.
Counting the numbers of all
possible $\sigma_{1}$ (hence $k$ and $m_{*}=(\mu_{*},j_{*},k_{*})$)
and the polynomials $p\in\mathcal{P}_{K}^{(1)}$,
the total excision measure for $\omega$ is
\[
\lessdot~ N^{d} K^{d} K^{5+\frac{6C_{1}}{C_{3}}} \Phi(K)^{-\frac{1}{2}}\eta
= N^{-(\frac{C_{1}}{2C_{5}}-d-\frac{d+5+\frac{6C_{1}}{C_{3}}}{C_{5}})}\eta
<\frac{\eta}{N^{4}},
\]
in view of \eqref{cons 2}, or equivalently,
\begin{equation*}
  \frac{C_{1}}{2C_{5}}-d-\frac{d+5+\frac{6C_{1}}{C_{3}}}{C_{5}}>4.
\end{equation*}
Define
\begin{equation*}
  \mathscr{V}_{j+1}^{(1)}= \bigcup_{ M^{j}< N\leq M^{j+1}}\bigcup_{|k|\leq N}\bigcup_{m_{*},
  ,|k_{*}|\leq K} \bigcup_{p\in\mathcal{P}^{(1)}_{K}}\left\{\omega: |p(\sigma_{1})|\leq \frac{\eta}{\sqrt{\Phi(K)}}\right\}.
\end{equation*}
Then we have
\begin{equation}\label{V 1}
  \textrm{mes}[\scr{V}_{j+1}^{(1)}]\leq \eta\sum_{M^{j} <N\leq M^{j+1}}
N^{-4}\leq \frac{M^{j+1}}{M^{4j}}< \frac{1}{100}\cdot\frac{\eta}{4^{j}}.
\end{equation}

\subsection{Separation property.}\label{IL 2}

Now we study the separation property
\eqref{sepa ind.} in the iteration lemma but also for $T=\scr{F}'[y_{j}]$.
Recall the definition of $T^{\sigma}$ in \eqref{T(sigma)}.
Suppose both $(T^{\sigma_{1}}_{N})^{-1}$ and
$(T^{\sigma_{2}}_{N})^{-1}$ fail \eqref{sepa ind.}.
From the analysis in our main construction section \ref{constr. II},
there should exist some
$m_{1}=(\mu_{1},j_{1},k_{1}), m_{2}=(\mu_{2},j_{2},k_{2})$
and $p_{1},p_{2}\in \mc{P}_{N}^{(1)}$
 such that
 \begin{equation*}
 	|p_{i}(\sigma_{i}')|<\frac{\eta}{\sqrt{\Phi(N)}},\quad
 	\sigma_{i}'=\sigma_{i}+\langle k_{i},\omega\rangle+\mu_{i}\lambda_{j_{i}},\quad i=1,2.
 \end{equation*}
 Recall that $p_{1}$ and $p_{2}$ are linear functions
 taking the form of
 \begin{equation*}
 	p_{i}(\sigma)=\sigma+a_{i}(\omega),\quad
 	|\partial_{\omega}a_{i}(\omega)|\leq C,\quad
 	i=1,2.
 \end{equation*}
Then we see that
\begin{equation*}
|p_{1}(\sigma_{1}')-	p_{2}(\sigma_{2}')|
=|\langle k+k_{1}-k_{2},\omega\rangle+\mu_{1}\lambda_{j_{1}}-
\mu_{2}\lambda_{j_{2}}+a_{1}(\omega)-a_{2}(\omega)|
< \frac{2\eta}{\sqrt{\Phi(N)}},
\end{equation*}
with
\begin{equation*}
2N<4N-2N<|k+k_{1}-k_{2}|	<2 N^{C_{3}}.
\end{equation*}

Take parameter excision for $\omega$ as follows.
Define by
$\mathcal{P}_{N}^{(2)}$  the set of the following polynomials
\begin{align*}
&\tilde{p}(\sigma)=\sigma+ a_{0,1}-a_{0,2} ,\\
& a_{0,i}=p_{i}(\sigma)-\sigma,\quad \textrm{for some}\quad p_{i}\in\mathcal{P}_{N}^{(1)},\quad i=1,2,
\end{align*}
and consequently $\# ~\mathcal{P}_{N}^{(2)}\leq (\#~ \mathcal{P}_{N}^{(1)})^{2}$.
We estimate the measure of
\begin{equation*}
\mathscr{V}_{j+1}^{(2)}=\bigcup_{M^{j}<N\leq M^{j+1}}\bigcup_{\mu_{i},j_{i}, i=1,2}\bigcup_{2N<|k|< 2N^{C_{3}}}\bigcup_{p\in\mc{P}^{(2)}_{N}}
\left\{\omega: |p(\langle k,\omega\rangle+\mu_{1}\lambda_{j_{1}}-\mu_{2}\lambda_{j_{2}})|<
\frac{2\eta}{\sqrt{\Phi(N)}}\right\}	,
\end{equation*}
which satisfies
\begin{equation}\label{V 2}
\textrm{mes}[\scr{V}_{j+1}^{(2)}]\lessdot
 \eta\sum_{M^{j}<N\leq M^{j+1}} N^{-(\frac{C_{1}}{2}-d C_{3}-10-\frac{12 C_{1}}{C_{3}})}
< \eta\sum_{M^{j}<N\leq M^{j+1}} N^{-4}<\frac{1}{100}\cdot\frac{\eta}{4^{j}},
\end{equation}
in view of
\eqref{cons 3}.

At step  $j+1$, we have excluded parameters twice.
One is to control $T_{N}^{-1}$ and
the other one is to establish the separation property. Denote
by $\mathscr{U}_{j+1}$ the "good" parameter set such that the
analysis in the $(j+1)$-th step is valid, i.e.,
\begin{equation*}
  \mathscr{U}_{j+1}=\scr{U}_{j}\setminus(\scr{V}_{j+1}^{(1)}\cup \scr{V}_{j+1}^{(2)}).
\end{equation*}
Then we obtain from \eqref{V 1} and \eqref{V 2} that
\begin{equation*}
\textrm{mes}[\mathscr{U}_{j}\setminus \mathscr{U}_{j+1}]< \frac{1}{50}\cdot
\frac{\eta}{4^{j}},
\end{equation*}
which verifies statement \textbf{(S4)}$_{j+1}$.

\subsection{Approximate solution and new error.}\label{IL 3}
Consider the $(j+1)$-th step and take $N=M^{j+1}$. Let
$$
y_{j+1}=y_{j}+\Delta_{j},
$$
where  $\Delta_{j}$ is given by the Newton equation
\begin{equation*}
\Delta_{j}=-(T_{N})^{-1}\Gamma_{10 M^{j}}\scr{F}[y_{j}],\quad T=\scr{F}'[y_{j}].
\end{equation*}
By \eqref{coup 1} and $M^{c}\approx 1$, we obtain
\begin{equation}\label{Delta 0}
  \|\Delta_{j}\|\leq \|(T_{N})	^{-1}\| \cdot \|\scr{F}[y_{j}]\|
\leq 2\Phi(\frac{1}{2} N) e^{-2 M^{jc}}< e^{-\frac{3}{2} N^{c}},
\end{equation}
and consequently
\begin{equation}\label{Delta 1}
  |\Delta_{j}(m)|< e^{-|k|^{c}} e^{-\frac{1}{2} N^{c}},
\end{equation}
since $|k|\leq N$.

Consider $\partial_{\omega}\Delta_{j}=-(\partial_{\omega}(T_{N})^{-1})\cdot \scr{F}[y_{j}]-(T_{N})^{-1} \partial_{\omega} \scr{F}[y_{j}]$.
The estimate for $(T_{N})^{-1} \partial_{\omega} \scr{F}[y_{j}]$
remains the same to $\Delta_{j}$ and thus
we have
\begin{equation}\label{Delta 2}
	|(T_{N})^{-1} \partial_{\omega} \scr{F}[y_{j}](m)|
	< e^{-|k|^{c}} e^{-\frac{1}{2} M^{(j+1)c}}.
\end{equation}
Note that
\begin{equation*}
\partial_{\omega} (T_{N})^{-1}= -(T_{N})^{-1}(\partial_{\omega}
T_{N}) (T_{N})^{-1}	
\end{equation*}
and similar arguments in \eqref{T'_omega} yield
\begin{equation*}
\|(\partial_{\omega} T_{N})\|\lessdot  N^{2d}.
\end{equation*}
One readily finds that the polynomial growth (in $N$) of
$\partial_{\omega}(T_{N})^{-1}$ can always be controlled
by the exponential decay of $\|\scr{F}[y_{j}]\|$. Then, by
shrinking the "analytic" strip, there is also
\begin{equation}\label{Delta 3}
	|\partial_{\omega}\Delta_{j}(m)|< e^{-|k|^{c}} e^{-\frac{1}{2} M^{(j+1)c}}.
\end{equation}

Combining \eqref{Delta 1}-\eqref{Delta 3}, we have
\begin{equation*}
  |\partial^{\alpha}_{\omega} \Delta_{j}(m)|< e^{-|k|^{c}} e^{-\frac{1}{10} M^{(j+1)c}},\quad |\alpha|=0,1,
\end{equation*}
and consequently
\begin{equation*}
  |\partial^{\alpha}_{\omega} y_{j+1}(m)|\leq e^{-|k|^{c}}\sum_{i=j_{0}}^{j}
e^{-\frac{1}{10} M^{c(i+1)}}
< e^{-|k|^{c}}.
\end{equation*}
Moreover, with $j_{0}$ large enough, we can also ensure that $\|y_{j+1}\|$ stays in a
small neighborhood of zero.
This verifies statement \textbf{(S1)}$_{j+1}$.

Newt we turn to estimate the new error.
We split $\scr{F}[y_{j+1}]$ into several terms
\begin{equation}\label{error 0}
\begin{aligned}
\scr{F}[y_{j+1}]=&\scr{F}[y_{j}+\Delta_{j}]=\scr{F}[y_{j}]
+D\scr{F}[y_{j}]\Delta_{j}+R	\\
=&\scr{F}[y_{j}]+ T_{N}\Delta_{j}
+ (T-T_{N})\Delta_{j}
+ R\\
=& \scr{F}[y_{j}]-\Gamma_{10M^{j}}\scr{F}[y_{j}]\\
&+ \Gamma_{10 M^{j}}\scr{F}[y_{j}]+T_{N}\Delta_{j}\\
&+ (T-T_{N}) \Delta_{j}\\
&+ R,
\end{aligned}
\end{equation}
where
$$
R=\int_{0}^{1}\int_{0}^{1} \scr{F}''[y_{j}+st \Delta_{j}]
t \Delta_{j}^{\otimes 2}\textrm{d}s\textrm{d}t.
$$

Observe that $\Gamma_{10 M^{j}} \scr{F}[y_{j}]+T_{N}\Delta_{j}=0$ and
by \eqref{Delta 0}
\begin{equation}\label{error 1}
	\|R\|\lessdot\|\Delta_{j}\|^{2}
	<\frac{1}{4} e^{-2 M^{(j+1)c}}.
\end{equation}

We further decompose $(T-T_{N})\Delta_{j}$ into
\begin{align*}
(T-T_{N}) \Delta_{j}=&(1-\Gamma_{N})T \Delta_{j}
= (1-\Gamma_{N})T~ \Gamma_{\frac{N}{2}}\Delta_{j}
+(1-\Gamma_{N})T~ (1-\Gamma_{\frac{N}{2}})\Delta_{j}\\
=& (1-\Gamma_{N})T~ \Gamma_{\frac{N}{2}}\Delta_{j}
+ (1-\Gamma_{N})T~ (1-\Gamma_{\frac{N}{2}}) (T_{N})^{-1}
\Gamma_{10 M^{j}} \scr{F}[y_{j}] .
\end{align*}
Using the off decay estimate
\begin{equation*}
|T(m,m')|< \epsilon e^{-|k-k'|^{c}},\quad k\neq k'
\end{equation*}
for $T=\scr{F}'[y_{j}]$, we have
\begin{equation*}
\begin{aligned}
	\|(1-\Gamma_{N})T~\Gamma_{\frac{N}{2}}\|=&
	\sup_{\|y\|=1}\|(1-\Gamma_{N}) T~\Gamma_{\frac{N}{2}}
	y\|\\
	\leq & \sup_{\|y\|=1} \left(\sum_{|k|\geq N}(\sum_{|k'|\leq \frac{N}{2}}
	|T(m,m')|\cdot |y(m')|)^{2}\right)
	^{1/2}
	\leq  \frac{1}{8} e^{-\frac{1}{2} N^{c}}
	\end{aligned}
\end{equation*}
since $|k-k'|\geq \frac{N}{2}$.
Then there is
\begin{equation}\label{error 2}
	\|(1-\Gamma_{N})T~ \Gamma_{\frac{N}{2}}\Delta_{j}\|
	\leq \frac{1}{8} e^{-\frac{1}{2} N^{c}} e^{-\frac{3}{2}N^{c}}=\frac{1}{8} e^{-2 M^{(j+1)c}}.
\end{equation}

Similarly, we obtain
\begin{equation*}
	\|(1-\Gamma_{N})T~ (1-\Gamma_{\frac{N}{2}}) (T_{N})^{-1}
\Gamma_{10 M^{j}}\|\leq \frac{1}{8} e^{-\frac{1}{4} N^{c}}
\end{equation*}
and hence
\begin{equation}\label{error 3}
	\|(1-\Gamma_{N})T~ (1-\Gamma_{\frac{N}{2}}) (T_{N})^{-1}
\Gamma_{10 M^{j}} \scr{F}[y_{j}]\|\leq
\frac{1}{8} e^{-(\frac{1}{4}+\frac{2}{M^{c}}) N^{c}}
\leq \frac{1}{8} e^{-2 M^{(j+1)c}}.
\end{equation}

Finally, we show the estimate for $(1-\Gamma_{10 M^{j}}) \scr{F}[y_{j}]$.
Since $\textrm{Supp}~y_{j}\subset\mc{L}_{M^{j}}$, we obtain
from $\scr{F}=D+\epsilon \scr{W}$ that
\begin{equation*}
	(1-\Gamma_{10 M^{j}})\scr{F}[y_{j}]=(1-\Gamma_{10 M^{j}})
	\scr{W}[y_{j}].
\end{equation*}
Furthermore, taking only the $k$-component into consideration, we have
\begin{align*}
\|(1-\Gamma_{10 M^{j}})\scr{W}[y_{j}]\|\lessdot \sum_{\zeta\geq 10}
\sum_{|\alpha|+|\beta|=\zeta}	|b_{\alpha,\beta}|
\sum_{\substack{|\sum_{1\leq s\leq \zeta}k_{s}|\geq 10M^{j}\\
|k_{\zeta}|\leq M_{j}, 1\leq s\leq \zeta}} \|y_{j}\|^{\zeta} e^{-\sum_{1\leq s\leq \zeta} |k_{s}|^{c}},
\end{align*}
where $b_{\alpha,\beta}$ corresponds to the coefficients in the power series expansion of $f$ in \eqref{schro-like}
\begin{equation*}
  f(x,\bar x)=\sum_{\zeta\geq 0} \sum_{|\alpha|+|\beta|=\zeta}	b_{\alpha,\beta}
  x_{1}^{\alpha_{1}}\cdots x_{n}^{\alpha_{n}} \bar{x}_{1}^{\beta_{1}}\cdots
  \bar{x}_{n}^{\beta_{n}}.
\end{equation*}

Since the infimum of $\sum_{1\leq s\leq \zeta}(|k_{s}|/M^{j})^{c}$ with constraints
$\sum_{1\leq s\leq \zeta} (|k_{s}|/ M^{j})\geq 10$,
$(|k_{s}|/M^{j})\leq 1$ and $\zeta\geq 10$
is greater than three, we have
\begin{equation*}
	e^{-\sum_{1\leq s\leq \zeta}|k_{s}|^{c}}\leq e^{-3 M^{jc}}
	\leq e^{-\frac{3}{M^{c}} M^{(j+1)c}}.
\end{equation*}
Then by the smallness of $\|y_{j}\|$ (staying inside the analyticity domain of the function $f$),
we have
\begin{equation}\label{error 4}
\|(1-\Gamma_{10 M^{j}})\scr{F}[y_{j}]\|	\leq \frac{1}{4}
e^{-2 M^{(j+1)c}}
\end{equation}
provided $M^{c}\approx 1$.

Combining the estimates \eqref{error 1}-\eqref{error 4} for the split \eqref{error 0},
we have
\begin{equation*}
\|\scr{F}[y_{j+1}]\|	\leq e^{-2 M^{(j+1)c}}.
\end{equation*}
The estimate of $\|\partial_{\omega}\scr{F}[y_{j+1}]\|$ follows the same line and
we do not carry it our here.

\subsection{Final reckoning}\label{IL 4}

To complete the inductions, it suffices to apply the Neumann series
to  establish induction statements \textbf{(S3)}$_{j+1}$ and
\textbf{(S4)}$_{j+1}$. Let
$T_{j+1}=\scr{F}'[y_{j+1}]$ and $T_{j}=\scr{F}'[y_{j}]$,
whose restrictions on $\mathcal{L}_{N}$ are denoted by
by $T_{j+1;N}$ and $T_{j;N}$ respectively.
Note that
\begin{equation*}
\|\scr{F}'[y_{j+1}]-\scr{F}'[y_{j}]\|\lessdot \|y_{j+1}-y_{j}\|
\leq e^{-\frac{3}{2} M^{(j+1) c}}.
\end{equation*}
and thus for any $N\leq M^{j+1}$,
\begin{equation*}
    (T_{j+1; N})^{-1}=\left(1+\sum_{s=1}^{\infty} [(T_{j;N})^{-1} (T_{j+1;N}-T_{j;N})]^{s}\right)
    (T_{j;N})^{-1},
\end{equation*}
which is bounded by
\begin{equation*}
  \|(T_{j+1; N})^{-1}\|\leq \left(1+\sum_{s\geq 1} e^{-\frac{3}{2} M^{(j+1)c}} 2\Phi(\frac{1}{2} M^{j+1}) \right)\cdot 2\Phi(\frac{N}{2})
  < \Phi(N).
\end{equation*}
Moreover,  we find that
\begin{equation*}
  \left|\left[(T_{j;N})^{-1}(T_{j+1;N}-T_{j;N})^{s}\right](m,m')\right|\leq
  \Phi(N) e^{-\frac{3 s}{2} N^{c}}<\frac{1}{2^{s}} e^{-\frac{3}{4} (2N)^{c}}
  < \frac{1}{2^{s}} e^{-\frac{3}{4} |k-k'|^{c}},\quad
  s\geq 1,
\end{equation*}
and consequently by \eqref{coup 2}
\begin{equation*}
  |(T_{j+1;N})^{-1}(m,m')|<e^{-\frac{1}{2}|k-k'|^{c} },\quad |k-k'|> N^{\frac{1}{C_{6}}}.
\end{equation*}

The separation property for $T_{j+1}^{\sigma}$ follows the same way by
applying the Neumann series and we do not repeat it here.

\bibliographystyle{plain}
\def\cprime{$'$} \def\cprime{$'$} \def\cprime{$'$}

\end{document}